\pdfoutput=1
\newif\ifreport
\newif\ifnotreport
\reporttrue % Comment for paper version.
\ifreport
  \notreportfalse
\else
  \notreporttrue
\fi

\newif\ifreview
\newif\ifnotreview

\ifreport
  \reviewtrue % Proofs are available in the report.
\else
\fi

\ifreview
  \notreviewfalse
\else
  \notreviewtrue
\fi

\ifreview
\newenvironment{myproof}{\begin{proof}}{\end{proof}}
\newenvironment{detachedproof}[1]{\noindent\hspace{2em}{\itshape Proof of #1:}}{\hspace*{\fill}~\QED\par\endtrivlist\unskip}
\else
\newenvironment{myproof}{\expandafter\comment}{\expandafter\endcomment}

\fi

\ifreport
	\documentclass[10pt,a4]{article}
	
	\usepackage[utf8]{inputenc}
	
	\usepackage{cite}
	
	\usepackage[left]{lineno}
	
	\usepackage{microtype}
	\DisableLigatures[f]{encoding = *, family = * }
	
\else
	\documentclass[letterpaper, 10 pt, conference]{ieeeconf}  % Comment this line out if you need a4paper

	\IEEEoverridecommandlockouts                              % This command is only needed if
\fi

\usepackage{booktabs} % For formal tables

\usepackage[utf8]{inputenc}

\usepackage[T1]{fontenc}
\usepackage{ae,aecompl}

\usepackage{graphicx}

\usepackage{wrapfig}

\usepackage[compatibility=false]{caption}
\usepackage{subcaption}

\usepackage[english]{babel}

\usepackage{rotating}

\usepackage{verbatim}
\usepackage{bera,listings}

\usepackage{moreverb,url}

\usepackage{amsfonts}
\usepackage{amsmath}
\usepackage{amssymb}
\usepackage{stmaryrd}
\usepackage{bm}
\ifreport
	\usepackage{amsthm}
\fi

\usepackage{algorithm2e}

\usepackage{xcolor}

\usepackage{tabularx}
\usepackage{multirow}

\RequirePackage{tkz-graph}
\RequirePackage{tikz}
\usetikzlibrary{arrows,matrix,decorations.pathreplacing,positioning,chains,fit,shapes,calc,3d}

\usepackage{paralist}

\ifreport
  \usepackage[pdftex,colorlinks=true,urlcolor=blue,citecolor=black,anchorcolor=black,linkcolor=black]{hyperref}
\fi

\usepackage[noabbrev, capitalise, nameinlink]{cleveref}
\crefname{algocf}{algorithm}{algorithms}
\Crefname{algocf}{Algorithm}{Algorithms}

\graphicspath{{imgs/}}

\newif\ifcomments

\newcommand{\claudio}[1]{%
\ifcomments %
\textbf{\textcolor{red}{\scriptsize Claudio: #1}} %
\else %
\fi%
}%

\newcommand{\benoit}[1]{%
\ifcomments %
\textbf{\textcolor{blue}{\scriptsize Benoit: #1}} %
\else %
\fi%
}%

\newcommand{\footurl}[1]{\footnote{\url{#1}}}

\newcommand{\jsr}[0]{\ensuremath{ \hat{\rho} }}

\newcommand{\brackets}[1]{\ensuremath{ \left[ #1 \right] }}
\newcommand{\tuple}[1]{\ensuremath{ \left\langle #1 \right\rangle }}
\newcommand{\set}[1]{\ensuremath{ \left\{ #1 \right\}}}

\newcommand{\pargroup}[1]{\ensuremath{ \left( #1 \right)}}

\newcommand{\dert}[1]{\ensuremath{ \dot{#1} }}

\newcommand{\setreal}[0]{\ensuremath{\mathbb{R}}}

\newcommand{\norm}[1]{\left\lVert#1\right\rVert}
\newcommand{\bnorm}[1]{\big\lVert#1\big\rVert}

\newcommand{\opt}[1]{{#1}^\star}
\newcommand{\vectorOne}[1]{\brackets{%
\begin{matrix}
  #1
 \end{matrix}%
}}
\newcommand{\vectorTwo}[2]{\brackets{%
\begin{matrix}
  #1 \\
  #2
 \end{matrix}%
}}

\newcommand{\vectorFour}[4]{\brackets{%
\begin{matrix}
  #1 \\
  #2 \\
  #3 \\
  #4
 \end{matrix}%
}}

\newenvironment{aligneq*}%
{
\begin{equation*}
\begin{aligned}
}{
\end{aligned}
\end{equation*}
}

\newenvironment{aligneq}%
{
\begin{equation}
\begin{aligned}
}{
\end{aligned}
\end{equation}
}

\newcommand{\sys}[1]{System~(\ref{#1})}

\DeclareMathOperator*{\argmax}{arg\,max}

\newcommand{\Mats}{\mathcal{A}}
\newcommand{\Autom}{\mathbf{G}}

\newcommand{\Paths}[1][k]{\Autom_{#1}}
\newcommand{\Cycles}[1][k]{\Paths[#1]^\circ}

\newcommand{\Lang}{\mathcal{L}}
\newcommand{\Deb}[1][k]{\Autom^{[#1]}}
\newcommand{\DLang}[1][k]{\Lang^{[#1]}}
\newcommand{\DPaths}[1][t]{\Autom_{#1}^{[k]}}

\newcommand{\LangAutom}[1][\Autom]{{#1}^*}
\newcommand{\Oracle}{\mathcal{O}}

\newcommand{\DBA}{De Bruijn graph}

\newcommand{\Tr}{\top}

\newtheorem{definition}{Definition}{}
\newtheorem{example}{Example}{}
\newtheorem{lemma}{Lemma}{}
\newtheorem{remark}{Remark}{}
{}
\newtheorem{proposition}{Proposition}{}
{}
\newtheorem{problem}{Problem}{}
\newtheorem{theorem}{Theorem}{}
{}
{}

\begin{document}

\ifnotreport
	\bstctlcite{rmurl}
\fi

\newcommand{\Title}{Minimally Constrained Stable Switched Systems \\ and Application to Co-simulation}
\newcommand{\AuthA}{Cláudio Gomes}
\newcommand{\AuthB}{Raphaël M. Jungers}
\newcommand{\AuthC}{Benoît Legat}
\newcommand{\AuthD}{Hans Vangheluwe}

\ifreport
	
	\vspace*{0.35in}

	% title goes here:
	\begin{flushleft}
	
		{\Large
		\textbf\newline{\Title}
		}
		\newline
		% authors go here:
		\\
		\begin{center}
			\AuthA\textsuperscript{1,3},
			\AuthB\textsuperscript{2}, 
			\AuthC\textsuperscript{2}, \\
			\AuthD\textsuperscript{1,3,4}
		\end{center}
		
		\bigskip
		\bf{1} University of Antwerp, Belgium
		\\
		\bf{2} UCLouvain, Belgium
		\\
		\bf{3} Flanders Make vzw, Belgium
		\\
		\bf{4} McGill University, Canada
		\\
		\bigskip
		* claudio.gomes@uantwerp.be
		
	\end{flushleft}
	
\else
	\title{\LARGE \bf \Title}

	\author{Cláudio Gomes$^{1}$, Raphaël M. Jungers$^{2}$, Benoît Legat$^{3}$, and Hans Vangheluwe$^{4}$%
			\thanks{$^{1}$C. G. is a FWO Research Fellow, at the University of Antwerp, supported by the Research Foundation - Flanders (File Number 1S06316N). {\tt\small claudio.gomes@uantwerp.be}}%
            \thanks{$^{2}$R. M. J. is a F.R.S.-FNRS Research Associate, supported by the French Community of Belgium, the Walloon Region and the Innoviris Foundation.}%
			\thanks{$^{3}$B. L. is a F.R.S.-FNRS Research Fellow, at Université Catholique de Louvain. }%
			\thanks{$^{4}$H. V. is a Professor at the University of Antwerp. His work is partially supported by Flanders Make vzw, the strategic research centre for the manufacturing industry.}%
	}
	\maketitle
	\thispagestyle{empty}
	\pagestyle{empty}
\fi

\begin{abstract}
  We propose an algorithm to restrict the switching signals of a constrained switched system in order to guarantee its stability,
  while at the same time attempting to keep the largest possible set of allowed switching signals.
	
  Our work is motivated by applications to (co-)simulation %and consensus algorithms,
  where numerical stability is a hard constraint,
  but should be attained by restricting as little as possible the allowed behaviours of the simulators.

  We apply our results to certify the stability of an adaptive co-simulation orchestration algorithm,
  which selects the optimal switching signal at run-time, as a function of (varying) performance and accuracy requirements. %, and we show that it is possible to compute the optimal stabilization for a class of consensus algorithms.
\end{abstract}

\section{Introduction}
\label{sec:introduction}

A switched system is defined as
\begin{aligneq}\label{eq:switched_system}
x_{k+1} &= A_{\sigma_k} x_{k} : \ \ \sigma_k \in \set{1, \ldots, m},\  A_{\sigma_k} \in \Mats
\end{aligneq}
where
$\set{1, \ldots, m}$ is the set of modes,
$\sigma_k$ is the mode active at time $k$, and
$\Mats = \set{A_1, \ldots, A_m} \subseteq \mathbb{R}^{n \times n}$ is a set of real matrices.
We denote the sequence $\sigma_0, \sigma_1, \ldots, \sigma_k$ as the \emph{switching signal}.
\ifreport
	$A_{\sigma_k} \in \Mats$ represents the matrix used to compute $x_{k+1}$ from $x_{k}$ at time $k$.
\fi

\ifreport
  Switched systems are widely used to model many dynamical systems in modern engineering
  including viral mutations in a patient's body \cite{hernandez2011discrete},
  \emph{trackability} of malicious agents in a sensor network \cite{jungers2008efficient},
  or scheduling of thermostatically controlled loads (TCLs) \cite{nilsson2017class}.
  
\else
  Switched systems are used to model many systems \cite{jungers2008efficient}.
\fi
In this paper, we are motivated by a new application \cite{Gomes2017d} in the field of co-simulation.
It is a numerical technique to couple multiple simulators, each simulating a part of a coupled system, in order to compute the overall behavior more efficiently %
\ifreport
	\cite{Kubler2000,Busch2012,Busch2012a,Gu2001a,Gu2004,Gomes2018}.
\else
	\cite{Gomes2018}.
\fi
\ifreport
	One of the objectives is to leverage different mature simulation tools, even when the details of each simulation tool are unknown \cite{Blochwitz2012,Blochwitz2011}, and it has been applied in fields such as automotive, electricity distribution, maritime, railway, etc\ldots
	See \cite{Schweizer2016,Schweizer2015c,Gomes2016,Gomes2017} and references thereof, for applications.
\fi

As illustrated in \cref{fig:cosim_overview}, in order to run a co-simulation, each simulator approximates the solution of a differential equation, exchanging values with other simulators at agreed-upon communication times.
Between communication time instants, the unknown inputs are approximated with extrapolation functions.

\ifreport
  \begin{figure}[htb]
  \begin{center}
    \includegraphics[width=0.5\textwidth]{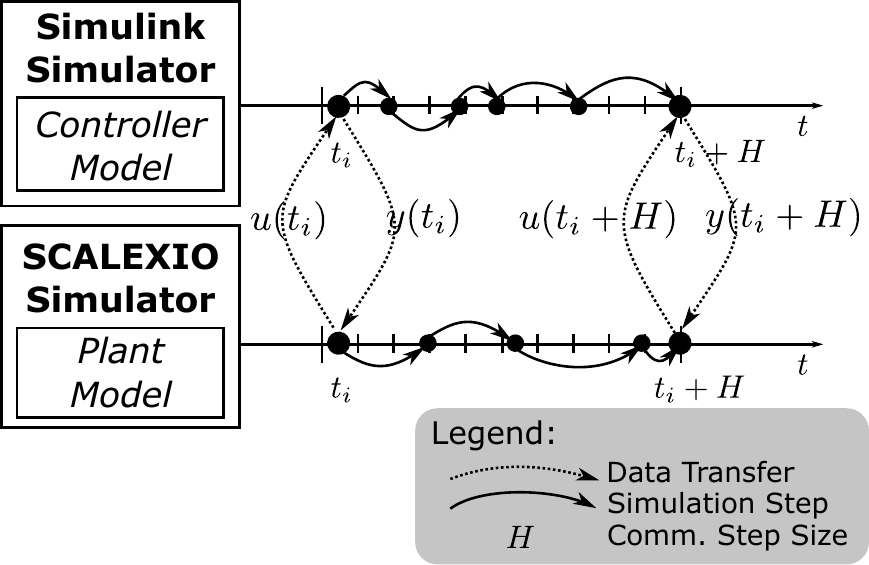}
    \caption{Co-simulation coordination example.}
    \label{fig:cosim_overview}
  \end{center}
  \end{figure}
\else
  \begin{wrapfigure}[10]{r}{0.3\textwidth}
  \begin{center}
    \includegraphics[width=0.3\textwidth]{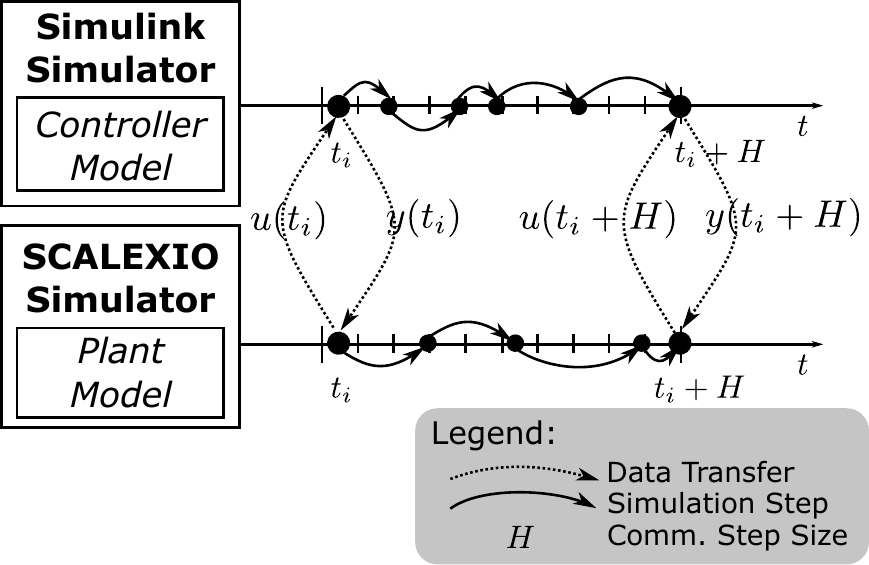}
    \caption{Co-simulation.}
    \label{fig:cosim_overview}
  \end{center}
  \end{wrapfigure}
\fi

\ifreport
  To improve the performance, during a co-simulation each simulator can adapt its policy, by varying:
  \begin{compactenum}
  \item the discretization step size and/or the numerical algorithm (e.g., midpoint method, Runge-Kutta) as a result of error estimates;
  \item the input approximation function (e.g., by varying the Lagrange polynomial degree) as a result of the dynamics of the inputs; and
  \item the values exchanged and the order in which they are exchanged, as a result of the varying structure of the coupled system being simulated.
  \end{compactenum}
\else
  To improve performance, during a co-simulation each simulator can adapt its policy:
  \begin{inparaenum}[(i)]
  \item the discretization step size and/or the numerical algorithm 
  as a result of error estimates;
  \item the input approximation function 
  as a result of the dynamics of the inputs; and
  \item the order and the values exchanged, as a result of varying the structure of the coupled system being simulated.
  \end{inparaenum}
\fi
We will call \emph{policy sequence} to the sequence of policies taken by all simulators over time.

\ifreport
	The ability to adapt allows one to find the best tradeoff between accuracy and computation power, to meet the available time. The applications are not restricted to co-simulation.
	For example, in Model Predictive Control~\cite{Garcia1989}, the controller needs to be able to simulate the system in real-time.
\else
	For other applications of adaptive methods, see \claudio{cite report. but I'm not sure we actually describe this in the techreport. Maybe just give the model predictive control example.}
\fi

In the context of (co-)~simulation, it is fundamental to ensure that a (co-)~simulation method preserves the stability of the system being (co-)~simulated.
The error dynamics of a (co-)~simulation can be modelled as a switched system (\cref{sec:application}),
\ifreport
  hence the preservation of stability becomes a problem of deciding the stability of a switched system (\sys{eq:switched_system}).
\else
  hence this question becomes one of deciding stability of \sys{eq:switched_system}.
\fi
Furthermore, as the choice of future policies is influenced by the past policies, we consider
\emph{constrained switched systems}%
\ifreport
	, a recently developed framework allowing us to model the memory of the system (see \cref{sec:problem_formulation}).
\else
	.
\fi

If there exist one or more policy sequences that can make the (co-)~simulation method unstable, then the simulators have to be forbidden from following these.
In the context of constrained switched systems, there are many ways of forbidding a policy sequence, and each way also forbids sequences that do not make the (co-)~simulation unstable.

We tackle the problem of how to best forbid policy sequences that make the constrained switched system unstable.
We propose that the best solution is to maximize the \emph{entropy} of the stabilized constrained switched system in order to maximize the adaptability of the resulting (co-)~simulation method.

In the next section, we formulate the problem%
\ifreport
  \ of how to forbid ``bad'' policy sequences (that cause the system to be unstable), while minimizing the number of good sequences that are forbidden as side effect.
\else
	.
\fi
Then, in \cref{sec:stabilization_ss} we propose an algorithm that approximates the solution, and we prove that the it terminates and that the resulting system is stable. Furthermore, we provide a lifting technique that yields better solutions.
\ifreport
	\Cref{sec:application} explores applications of our contribution to simulation and co-simulation. \claudio{Should we include the consensus algorithms? I have omitted it, and it might be useful for the journal paper.}
\else
    \Cref{sec:application} describes applications.
\fi
\Cref{sec:related_work} presents related work and \cref{sec:conclusion} concludes.

\ifnotreport
  An extended version of this work is presented in \claudio{[Rep]}.
\fi

\section{Problem Formulation}
\label{sec:problem_formulation}

\ifreport
	We introduce the dynamical systems being considered, called \emph{Constrained Switched Systems} \cite{Philippe2016}, and we formulate our research problem.
\else
	In this section, we formulate the research problem, and introduce the necessary concepts.
\fi

In practice, some switching signals of \sys{eq:switched_system} may not be relevant, and a way to represent the sensible ones is required.
\ifreport

	For instance, in the switched system described in \cref{ex:adaptive_numerical_method}, it may not make sense that the Runge-Kutta method is used right after a Forward Euler. 
	This is because the convergence rate of each method is too different to warrant a switch without first taking a step with the Midpoint method.
	
	To model these constraints, we introduce the notion of constrained switched system.
	When compared to \sys{eq:switched_system}, constrained switched systems incorporate a representation of the allowed switching signals using an automaton \cite{Linz2011}.
\fi

\begin{definition}
  \label[definition]{def:caut}
  Given a bounded set of matrices $\Mats = \set{A_1, \ldots, A_m}$, we define an \emph{automaton} as a directed and labelled graph
  $\Autom = (V, E)$, with nodes $V$ and edges $E$
  such that no node has zero ingoing or outgoing degree.
  Each edge $(v, w, \sigma) \in E$ represents a transition from node $v \in V$ to node $w \in V$,
  where $\sigma \in \set{1, \ldots, m}$ is the \emph{label}, corresponding to $A_\sigma$.
\end{definition}
\ifreport
	An example automaton, illustrating possible constraints on the system described in \cref{ex:adaptive_numerical_method}, is shown in \cref{fig:order_adaptive_method}.
\fi

We say that the switching signal, or word, $s = \sigma_0 \sigma_1 \ldots \sigma_{k-1}$ is \emph{accepted} by an automaton $\Autom$ if it corresponds to a path in $\Autom$, that is, if there exists $v_0,v_1,\ldots,v_k \in V$, such that
$(v_j, v_{j+1}, \sigma_j) \in E$ for all $j = 0, \ldots, k-1$.
An accepted word induces an accepted matrix product $A_s = A_{\sigma_{k-1}} \cdots A_{\sigma_1} A_{\sigma_0} \in \Mats^k$.
\ifreport
	% For any $k>0$, the word $s^k$ is the concatenation of $s$ with itself $k-1$ times.
	% \claudio{The previous sentence is only needed for the consensus.}
	
	For example, the word $(\mathit{fe}, 0.001),(\mathit{fe}, 0.002),(\mathit{md}, 0.002)$ is accepted by the automaton shown in \cref{fig:order_adaptive_method}. 
	This word induces the matrix product $\tilde{A}_{\mathit{md}, 0.002} \tilde{A}_{\mathit{fe}, 0.002}  \tilde{A}_{\mathit{fe}, 0.001}$.
	
\fi

We denote the set of accepted words of length $k$ as $\Paths[k]$,
and the set of all words accepted by the automaton as $\LangAutom = \bigcup_{k=1}^\infty \Autom_k$.
Moreover, $\Cycles[k]$ denotes the set of accepted cycles of length $k$.

\ifreport
	
	For example, 
	\begin{aligneq*}
	 (\mathit{fe}, 0.001),(\mathit{fe}, 0.002),(\mathit{md}, 0.002) &\in \Paths[3] \text{, and} \\
	 (\mathit{fe}, 0.002),(\mathit{fe}, 0.001) &\in \Cycles[2].
  \end{aligneq*}
	
	One can see that given a word $\sigma(0) \ldots \sigma(k-1) \in \Paths[k]$, any sub-word
	$\sigma(i)\ldots\sigma(j)$ for any $0 \leq i \leq j < k$, satisfies
	$\sigma(i)\ldots\sigma(j) \in \Paths[j-i+1]$.
	Moreover, since every node has at least one outgoing edge in \cref{def:caut},
	for any $k' > k$, there exists $\sigma(k)\ldots\sigma(k'-1)$ such that $\sigma(0)\ldots\sigma(k'-1) \in \Paths[k']$.
\fi

\ifreport
	\begin{figure}[tbh]
	\begin{center}
    \includegraphics[width=0.8 \textwidth]{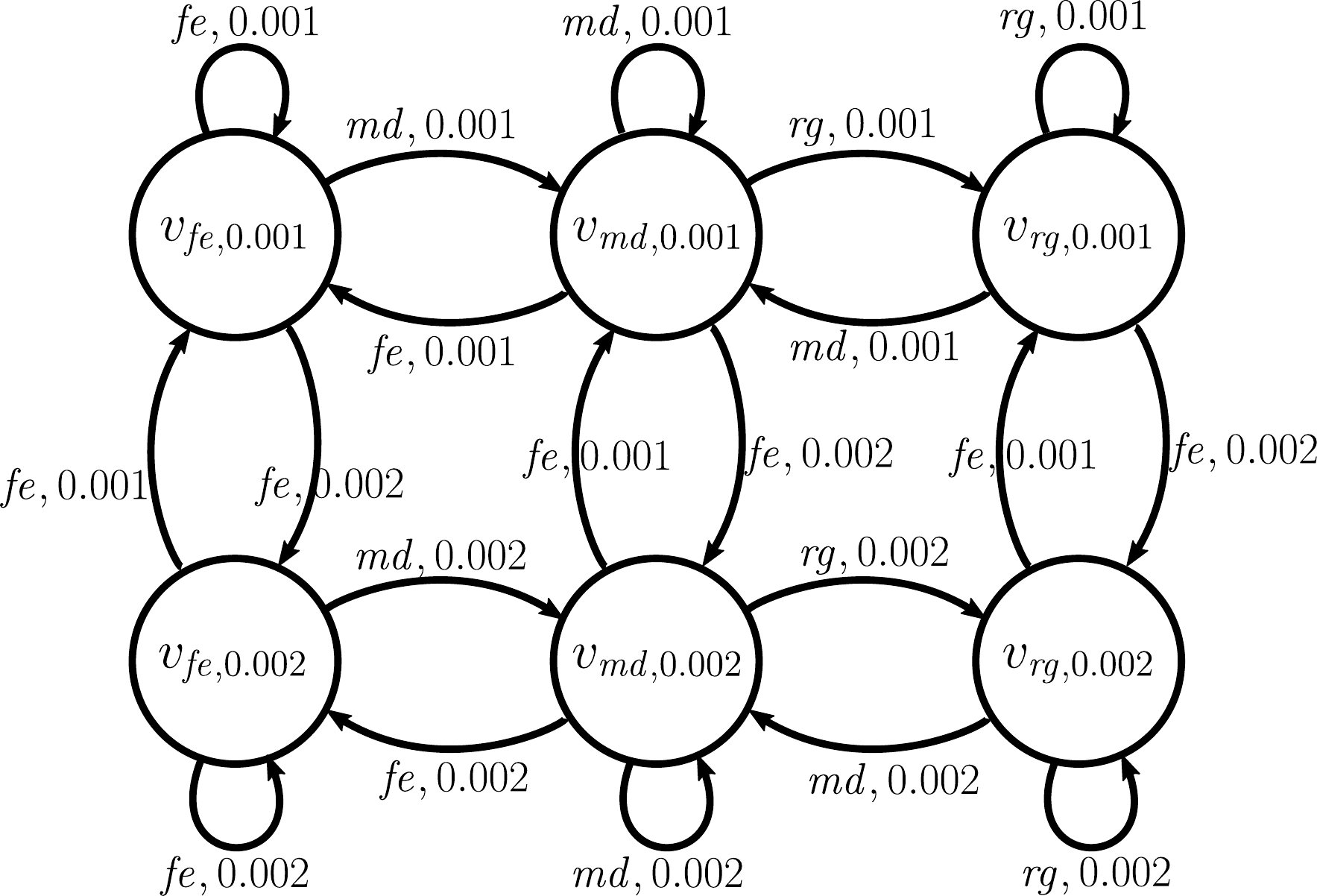}
    \caption{Example automaton for \cref{ex:adaptive_numerical_method}.}
    \label{fig:order_adaptive_method}
	\end{center}
	\end{figure}
\fi

\begin{definition}[CSS]
  Given a set of matrices $\Mats = \set{A_1, \ldots, A_m}$,
  and an automaton $\Autom = (V, E)$,
  we define a \emph{constrained switched system} (CSS) $S = \tuple{\Mats, \Autom}$
  as a system where the variable $x_k$ satisfies:
  \begin{equation}\label{eq:css}
    x_{k+1} = A_{\sigma_k} x_k: \quad \sigma_0 \ldots \sigma_{k-1} \in \Paths[k].
  \end{equation}
\end{definition}

We say that \sys{eq:css} is \emph{stable} iff
\begin{equation}\notag
    \lim_{k \to \infty} \norm{x_{k}} = \lim_{k \to \infty} \norm{ A_{\sigma_{k-1}} \cdots A_{\sigma_0} x_{0}} = 0,
\end{equation}
for any word $\sigma_0 \ldots \sigma_{k-1} \in \Paths[k]$ and any $x_0 \in \mathbb{R}^n$.

To determine the stability of a CSS, we introduce the \emph{constrained joint spectral radius}.
\begin{definition}[{\cite[Definition 1.2]{Dai2012}}]
The \emph{constrained joint spectral radius} is defined as
\begin{equation}\notag%\label{eq:cjsr}
\jsr(S) = \lim_{k \to \infty} \jsr_k(S)
\text{ where } \jsr_k(S) = \sup_{w \in \Paths[k]} \norm{A_w}^\frac{1}{k},
\end{equation}
and $\norm{ \cdot }$ is any matrix norm that satisfies the sub-multiplicative property.
\end{definition}

\begin{proposition}[{\cite[Lemma 3.1]{Dai2012}}]\label{prop:stability}
If $\jsr(S) < 1$, then the CSS is stable.
\end{proposition}

\ifreport
  It can be shown~\cite[Lemma~1.2]{Jungers2009} that
  $\lim_{k \to \infty} \jsr_k(S) = \inf_{k \geq 1} \jsr_k(S)$.
  Therefore, for any $k > 0$, $\jsr_k(S)$ is an upper bound to $\jsr(S)$.
  This fact, together with \cref{prop:stability}, gives us a way to check whether a given CSS $S$ is stable: 
  \begin{compactenum}
  \item Pick a finite $k>0$, and compute $\jsr_k(S)$;
  \item If $\jsr_k(S) < 1$, then $\jsr(S) < 1$ and $S$ is stable;
  \item Otherwise, pick a larger $k$ and try again.
  \end{compactenum}
  
  If the CSS is unstable, then the above procedure will never terminate.
\fi

\ifreport
  A way to prove that a CSS is unstable is to find a switching signal that causes the system to be unstable. 
  For example, by finding a cycle $c \in \Cycles$ with $\rho(A_c) \geq 1$, where $\rho(A_c)$ denotes the spectral radius of the matrix $A_c$ induced by the cycle.
  In other words, finding a matrix product that, when repeated forever, causes the system to be unstable.
  
  Unstable cycles can be found by brute force or branch-and-bound variants~\cite{gripenberg1996computing, guglielmi2008algorithm, jungers2014lifted}.
  Naturally, these methods look first for unstable cycles with a small length.
  However, finding longer cycles becomes prohibitively high (see \cref{sec:application}).
  
  For example, a simple (and naive) procedure to find a cycle is to pick a finite $k$, enumerate all cycles $c \in \Cycles[k]$, check whether $\rho(A_c) \geq 1$, and stop when one such cycle is found.
  
  The method introduced in \cite{Legat2017} works well for finding long cycles.
  To certify the stability of a given CSS, it solves a semidefinite program to compute polynomial Lyapunov functions of degree $2d$.
  If the program is infeasible, it uses the dual certificate of infeasibity to generate an infinite switching signal of guaranteed growth rate.
  Subwords of this signal can be used to find unstable cycles.
  As cycles are found along an infinite switching signal, finding \emph{long} unstable cycles is not particularly more difficult.
  Moreover, if no unstable cycle can be found, one can retry with polynomial Lyapunov functions of degree $2(d+1)$.
  There is an improved guarantee on the growth rate of the infinite switching signal as the degree increases.
  
  The following definition formalizes the spectral radius of cycle induced matrix products.
  
  \begin{definition}[{\cite[Definition 1.2]{Dai2012}}]
  The \emph{generalized spectral radius} of a CSS $S$ is defined as:
  \begin{equation}\label{eq:gsr}
    \rho(S) = \limsup_{k \to \infty} \rho_k(S)
    \text{ where }\rho_k(S) = \sup_{c \in \Cycles[k]} \rho(A_c)^\frac{1}{k} \\
  \end{equation}
  \end{definition}

  It follows \cite[Proposition 1.6]{Jungers2009} that, for finite $k>0$,
  $$
    \rho_k(S) \leq \rho(S) \leq \jsr(S) \leq \jsr_k(S).
  $$
  Moreover, since $\Mats$ is bounded, it is shown in \cite[Theorem A]{Dai2012} that $\rho(S) = \jsr(S)$.
  
  \begin{remark}\label{rem:undecidability}
    The above discussion about proving that a CSS is unstable focused on finding finite cycles, as opposed to infinite paths.
    In fact, there is no guarantee that if a CSS satisfies $\jsr(S) \geq 1$ (i.e., is unstable), then a cycle $c$ with finite length exists, with $\rho(A_c) \geq 1$ (see \cite[Section 2.4]{Jungers2009} and \cite[Theorem~2]{Blondel2000}).
    However, the systems we experimented with, either satisfy $\jsr(S) > 1$, or $\jsr(S) < 1$.
    For these, the following result was used.
  \end{remark}
  
\else
  A way to prove that a CSS admits an unstable switching signal is to find a $c \in \Cycles$ with $\rho(A_c) \geq 1$.
\fi

\begin{proposition}[{\cite[Theorem~2.3]{Jungers2009}}]
  \label{prop:instability}
  If $\jsr(S) > 1$, then there exists a cycle $c \in \Cycles$ of length $k$ that satisfies $\rho(A_c) \ge 1$.
\end{proposition}

\ifreport
	Our goal is to optimally modify a given CSS, by forbidding unstable switching signal cycles from the language it generates.
\fi
The problem of finding such cycles is outside the scope of our work (see \cite{Legat2016} for the algorithm we used%
\ifreport
	, and references thereof for algorithms with the same goal%
\fi
).
As such, we introduce the following definition, which represents any algorithm available for this purpose.
\begin{definition}[Oracle]\label{def:oracle}
  Given $\epsilon > 0$, we define a stability oracle $\Oracle_\epsilon: S \to \{\texttt{Stable}\} \cup \bigcup_{k=1}^\infty \Cycles$, where $S$ is a CSS.
  The oracle $\Oracle_\epsilon$ returns either \texttt{Stable} certifying that $\jsr(S) < 1$ or a cycle $c \in \Cycles$ such that $\rho(A_c)^{1/k} > 1 - \epsilon$.
\end{definition}

We emphasize that the oracle has a (slightly) imperfect behaviour: in case $1-\epsilon < \jsr(S) < 1$,
one cannot guarantee what the outcome of the oracle will be.
This imperfection is intentional (see \cref{rem:undecidability}), as it models the state of the art \cite{parrilo2007approximation}.
\Cref{prop:instability} ensures that if $\jsr(S) > 1 - \epsilon$, there exists a $k$ and a cycle $c \in \Cycles$ such that $\rho(A_c)^{1/k} > 1 - \epsilon$.

We now proceed to define the set of possible different switching signals that are admissible.

\begin{definition}[Admissible Regular Language]
  \label{def:reg}
  We say that $\Lang=\LangAutom$ is the language \emph{recognized} by the automaton $\Autom$.
  A language is \emph{regular} if it is recognized by a finite automaton.
  A language $\Lang$ recognized by an automaton $\Autom$ is \emph{admissible} for $\Mats$ if the constrained switched system $S = \tuple{\Mats, \Autom}$ satisfies
  $\jsr(S) < 1$.
\end{definition}

Let $\Lang_0$ denote the language recognized by the automaton $\Autom_0$ of a given $S = \tuple{\Mats, \Autom_0}$.
Informally, our goal is to find the ``largest'' regular language $\opt{\Lang} \subseteq \Lang_0$ that is admissible.
\ifreport
	For this optimization problem to be well defined we need to find a metric for the objective.
	This metric should be in accordance to the fact that given $\Lang \subseteq \Lang'$, the objective should favor $\Lang'$.
\fi
A widely used notion to describe the size of a regular language is that of Entropy.

\begin{definition}[Entropy {\cite[Definition~4.1.1]{Lind1995}}]
    Given a regular language $\Lang$ recognized by an automaton $\Autom$,
    we define the \emph{entropy} as
    \ifreport
      $$h(\Lang) = \lim_{k \to \infty} \frac{1}{k} \log_2|\Paths|.$$
  		In the above, $|\Paths|$ represents the number of words of length $k$ accepted by the automaton $\Autom$.
    \else
       $h(\Lang) = \lim_{k \to \infty} \frac{1}{k} \log_2|\Paths|.$    
    \fi
\end{definition}
We denote the entropy of the language $\LangAutom$ recognized by an automaton $\Autom$ as $\LangAutom[h](\Autom)$.
\ifnotreport
	This quantity can be computed as shown in \claudio{cite report}.
\fi

If $\Lang \subseteq \Lang'$, then $\Paths \subseteq \Paths'$ for any $k$, and so $h(\Lang) \leq h(\Lang')$.
Our problem can now be formulated.
\begin{problem}\label[problem]{eq:prob}
  Given a CSS $\tuple{\Mats, \Autom_0}$,
  find the language $\opt{\Lang}$ solution of the following optimization problem:
  \begin{align}
\notag
    \opt{\Lang} = \sup_{\Lang \text{ regular}} h(\Lang) \text{ s.t.}\\
\notag
    \Lang & \subseteq \Lang_0,\\
    \label{eq:admissible}
    \Lang & \text{ is admissible for }\Mats.
  \end{align}
  where $\Lang_0$ is the language recognized by $\Autom_0$.
\end{problem}

\begin{remark}
  In \cref{eq:prob}, we restrict our attention to regular languages.
  While there are examples that highlight the benefit of using non-regular languages (see \cref{ex:pumping}),
  in practice, one needs an \emph{efficient} way of generating accepted switching signals.
  \ifreport
	For instance, during a co-simulation, at any step, the simulators need to compute as quickly as possible the set of policies that can be taken (see \cite[Section 4.4]{Gomes2017d} for how this can be done).
  \fi
  Automata allow the decision procedure to be fast, with little memory.
  In addition, as hinted in \cref{ex:1cycle}, regular languages may be constructed to approximate an admissible language with entropy arbitrarily close to the entropy of the optimal solution, even if that optimal solution is a non-regular language.
\end{remark}

\begin{example}
  \label{ex:pumping}
  Consider $\Mats=\set{A_1, A_2}$, with $A_1 = 2$ and $A_2 = \frac{1}{2}$, and $\Autom = (V, E)$, where
  $V=\set{v_1}$ and $E = \set{ (v_1, v_1, 1), (v_1, v_1, 2)}$.
  \ifreport
    That is, $\Autom$ has the form 
    \begin{equation*}
    \begin{tikzpicture}
      \Vertex[L={$v_1$}]{1}
      \draw[thick,->] (1) to [out=155,in=205,looseness=7] node [midway, fill=white] {$1$} (1);
      \draw[thick,->] (1) to [out=-25,in=25,looseness=7] node [midway, fill=white] {$2$} (1);
    \end{tikzpicture}
    \end{equation*}
  \fi
  The optimal solution $\opt{\Lang}$ of (relaxed) \cref{eq:prob} should include every word that has more $1$s than $2$s.
  As shown in \cite[Example~1.73]{Sipser2013}, no automaton can be built that accepts this language.
\end{example}

\begin{figure}
  \centering
  \ifreport
    \includegraphics[width=0.6\textwidth]{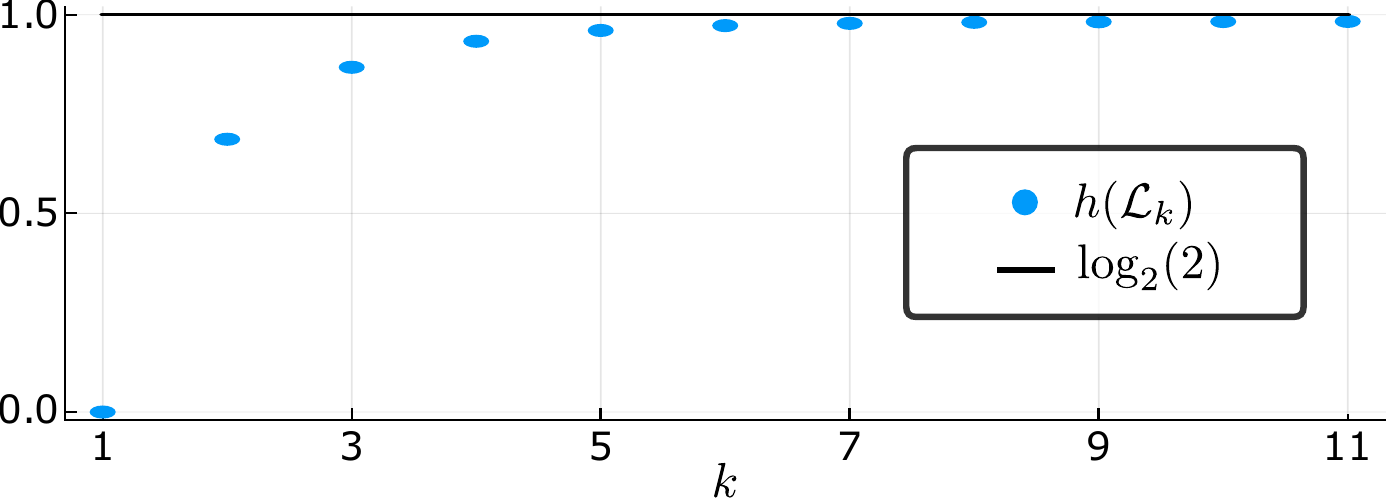}
  \else
    \includegraphics[width=0.3\textwidth]{1cycle_custom}  
  \fi
  \caption{Evolution of $h(\Lang_k)$ of \cref{ex:1cycle} in terms of $k$.}
  \label{fig:1cycle}
\end{figure}

\begin{example}
  \label{ex:1cycle}
  Consider $\Mats=\set{A_1, A_2}$, with $A_1 = 1$ and $A_2 = \frac{1}{2}$.
  A language is admissible if it does not contain the infinite repetition of the symbol 1.
  Let $\Lang_k$ be language of all words that do not contain $k$ consecutive 1's.
  \Cref{fig:1cycle} suggests that that $h(\Lang_k)$ tends to $\log_2(2)$ when $k$ tends to infinity.
  \ifreport
    The quantity $\log_2(2)$ denotes the entropy of the optimal solution.
  \fi
  % \benoit{Cite technical report where we prove in in the consensus section}
\end{example}

\section{Lift-and-Constrain Stabilization}
\label{sec:stabilization_ss}

\subsection{Constraining for more stability}
\label{sec:constrain_for_stability}

\Cref{alg:stabilization} details an iterative procedure that stabilizes a given CSS $S = \tuple{\Mats, \Autom}$, using the oracle in \cref{def:oracle}.
At each iteration, if the oracle returns a cycle $c=\sigma_k \ldots \sigma_k$, then $c$ is eliminated from $\Autom$.
The removal of a cycle can be accomplished by removing an edge of $\Autom$, thus potentially decreasing $\jsr(S)$.
After removing the cycle $c$, any infinite sequence in $\LangAutom$ for which $c$ is a subsequence will be eliminated too.
\ifreport
	This is illustrated in \cref{ex:edge_removal_example}.
\fi
The algorithm can produce an empty CSS, which does not imply that the original CSS is impossible to stabilize.
\ifreport
	An empty CSS is trivially stable.
\fi
\claudio{Maybe include an example?}

\begin{example}\label[example]{ex:edge_removal_example}
Consider the automaton in \cref{fig:edge_removal_example}, and suppose the oracle has returned the cycle $234$.
\ifreport
	This cycle is highlighted in red, in the figure.
\fi
Any of the edges in red can be removed to forbid the unstable sequence.
If edge $v_1 \xrightarrow{2} v_2$ is removed, the infinite sequences accepted by the resulting automaton end with either an infinite sequence of $2$'s, or an infinite sequence of $3$'s.
If edge $v_2 \xrightarrow{3} v_3$ is removed instead,
the resulting automaton accepts infinite sequences comprised of repeating subsequences which include $2$, or $3$, or $12$.
\end{example}

\begin{figure}[tbh]
\centering
  \begin{tikzpicture}[xscale=1.1,yscale=0.6]
    \Vertex[L={$v_1$}]{1}
        \EA[L=$v_2$,unit=3](1){2}
        \SOEA[L=$v_3$,unit=1.5](1){3}
        \SetUpEdge[style={post,{->}, thick},color=red]
        \Edge[label=$3$](2)(3)
        \Edge[label=$4$](3)(1)
        \SetUpEdge[style={post,bend right=18,thick},color=red]
        \Edge[label=$2$](1)(2)
        \SetUpEdge[style={post,bend right=18,thick}]
        \Edge[label=$1$](2)(1)

        \draw[thick,->] (2) to [out=-25,in=25,looseness=5] node [midway, fill=white] {$2$} (2);
        \draw[thick,->] (3) to [out=155,in=205,looseness=5] node [midway, fill=white] {$3$} (3);
  \end{tikzpicture}
  \caption{Automaton of Example~\ref{ex:edge_removal_example}.}
  \label{fig:edge_removal_example}
\end{figure}

As \cref{ex:edge_removal_example} shows, the choice of different edges to be removed has a different impact in the entropy of the resulting automaton.
\ifreport
	Informally, removing the edge $v_2 \xrightarrow{3} v_3$ seems to be the best choice because the resulting automata allows for \emph{more} sequences.
  This is corroborated by computing the entropy of the resulting automaton alternatives.
  See \cref{sec:computation_entropy} for how to compute the entropy in this example.
\fi

\begin{algorithm}
	 \KwData{A CSS $S=\tuple{\Mats, \Autom}$.}
	 \KwResult{A stable CSS $S=\tuple{\Mats, \Autom}$.}
     \While{$\Oracle_\epsilon(S) \neq {}$\texttt{Stable}}{
	 \begin{compactenum}
     \item\label{alg:stabilization:s2} Find $e \in \argmax \{\, \LangAutom[h](\Autom - e) \mid e \in E$, $e$ is\\an edge of the cycle $\Oracle_\epsilon(S) \,\}$\;
     \item Set $\Autom := \Autom - e$\;
	 \end{compactenum}
     }

   \caption{Stabilization algorithm for a constrained switched system. $\LangAutom[h](\Autom)$ denotes the entropy of the language recognized by $\Autom$. The difference $\Autom - e$ denotes the automaton obtained by removing the edge $e$ from $\Autom$.}
	 \label{alg:stabilization}
\end{algorithm}

\ifreport
	The following result demonstrates that \cref{alg:stabilization} always terminates.
\fi

\begin{theorem}\label{thm:termination}
  Given a CSS $S = \tuple{\Mats, \Autom}$ and
  an oracle satisfying \cref{def:oracle},
  \cref{alg:stabilization} terminates in finite time and the resulting CSS
  is stable.

  \begin{myproof}
    At each iteration of the algorithm, the number of edges of the automaton $\Autom = (V, E)$ decreases by one.
    Since at the beginning of the algorithm $|E|$ is finite,
    the algorithm must terminate after a finite number of iterations.
    The condition for termination of \cref{alg:stabilization} implies that the resulting system is stable.
  \end{myproof}
\end{theorem}

\ifreport
	\begin{remark}
	In \cref{thm:termination}, the assumption that the oracle in \cref{def:oracle} always terminates is crucial, as the problem solved by the oracle is undecidable in general (recall \cref{rem:undecidability}).
	\end{remark}
\fi

\subsection{Lifting for less conservativeness}

\ifreport
	\Cref{alg:stabilization} takes a constrained switched system $S=\tuple{\Mats, \Autom}$, and outputs a constrained switched system $S'=\tuple{\Mats', \Autom'}$ that is stable, while attempting to maximize the entropy of the language recognized by $\Autom'$.
	If we let $\Lang'$ denote this language, then, relating this to \cref{eq:prob}, $\Lang'$ is admissible and regular, and thereby a potential solution.
	However, it may not be the optimal solution.
	Similarly, if the algorithm returns an empty CSS, this does not mean that the original CSS is impossible to stabilize \claudio{Refer to the example, if we have one}.
\else
	The CSS returned by \cref{alg:stabilization} is not necessarily the optimal solution to \cref{eq:prob}.
\fi
To maximize the entropy of the stabilized CSS's, we propose to take an \emph{$M$-Path-Dependent lift} of the automaton representing the input language $\Lang_0$.

\begin{definition}[{\cite[Definition~3]{Philippe2016}}]\label[definition]{def:debruijn_graph}
  Given an automaton $\Autom$, we define the \emph{lifted} automaton $\Deb$ of degree $k$ as follows.
  For each path $v_0, \sigma_0, v_1, \sigma_1, \ldots, \sigma_{k}, v_{k+1}$ %
  \ifreport
    with length $k+1$
  \fi
   of $\Autom$,
  $\Deb$ has a node $u^- = v_0 \sigma_0 v_1 \sigma_1 \ldots \sigma_{k-1}v_k$,
  a node $u^+ = v_1 \sigma_1 v_2 \sigma_2 \ldots \sigma_{k}v_{k+1}$
  and an edge $(u^-, u^+, \sigma_k)$.
\end{definition}

\ifreport
	\Cref{fig:lifted_edge_removal_example} shows the second degree ($k=2$) lift of the automaton in \cref{fig:edge_removal_example}.
	
  \begin{figure}[tbh]
  \centering
    \includegraphics[width=0.6\textwidth]{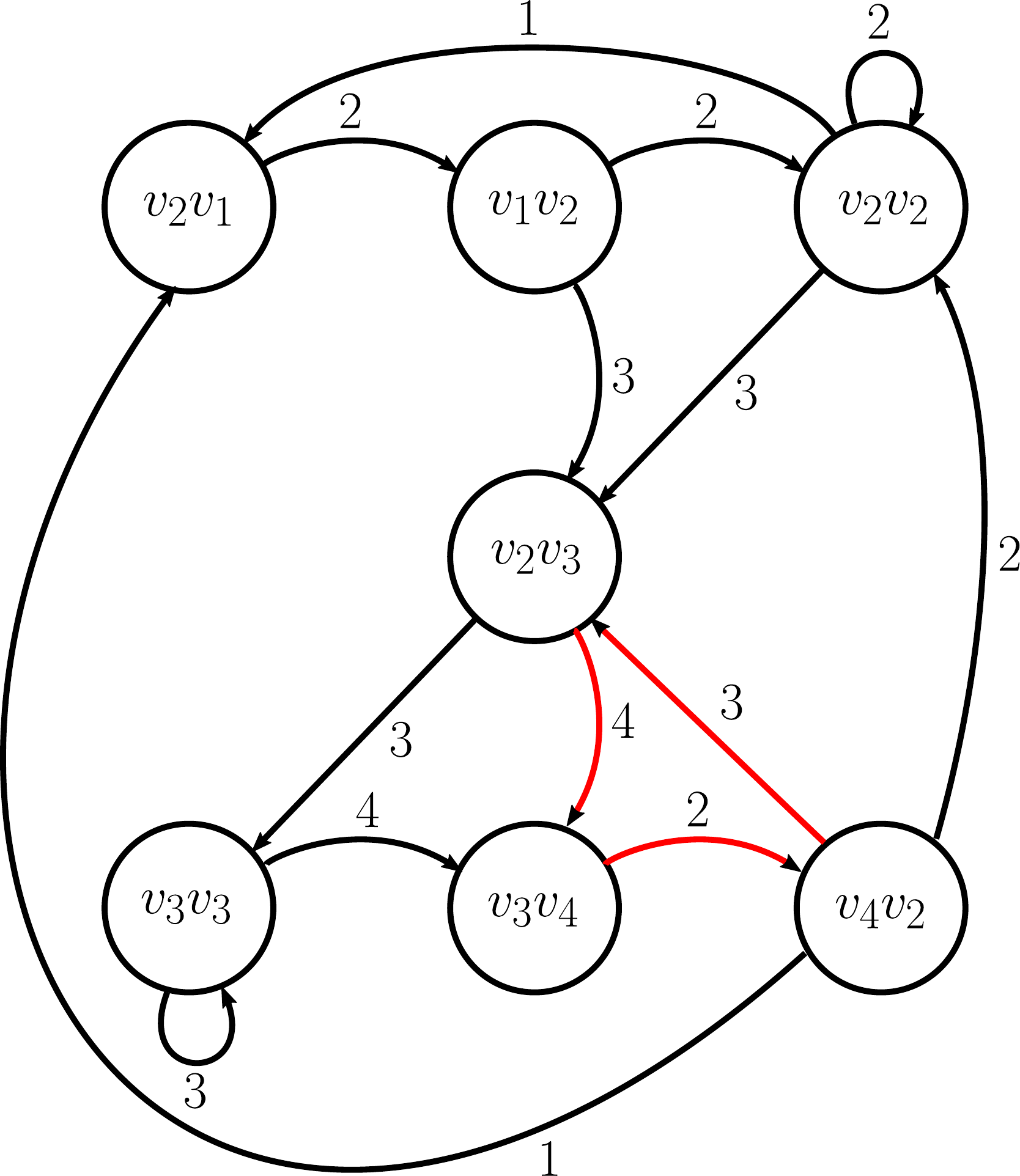}
    \caption{Second degree lifted automaton of Example~\ref{ex:edge_removal_example}.}
    \label{fig:lifted_edge_removal_example}
  \end{figure}
\fi

\ifreport
  The lifted automaton represents the same language, as shown by \cref{theo:incl}, but, as suggested by \cref{theo:remove_1_cycle} and illustrated by \cref{ex:1cycle}\claudio{Include a better example of this, using \cref{fig:lifted_edge_removal_example})}, lifting the automaton before applying \cref{alg:stabilization} allows one to obtain an admissible language with higher entropy.
\else
  The lifted automaton represents the same language, as shown by \cref{theo:incl}, but, as suggested by \cref{theo:remove_1_cycle} and illustrated by \cref{ex:1cycle}, lifting the automaton before applying \cref{alg:stabilization} allows to obtain an admissible language with higher entropy.
\fi

\begin{proposition}
  \label{theo:incl}
  Let $\Lang$ be the language recognized by $\Autom$ and $\DLang$ the language recognized by $\Deb$,
  where $\Deb$ is the lift of degree $k$ of $\Autom$.
  Then $\Lang = \DLang$.

  \begin{myproof}
    Consider a sequence $\sigma_{0} \ldots \sigma_{i-1}$.

    If $\sigma_{0} \ldots \sigma_{i-1} \in \Paths[i]$, there exists nodes $v_0, v_1, \ldots, v_i$ of $\Autom$ such that
    $$v_0, \sigma_0, v_1, \sigma_1, \ldots, \sigma_{i-1}, v_{i}$$ 
    is a path of $\Autom$.
    As no node has zero ingoing degree, there exists a path of length $k$ that ends in node $v_0$, denoted as $v_{-k}, \sigma_{-k}, \ldots, \sigma_{-1}, v_0$ in $\Autom$.
    By \cref{def:debruijn_graph}, for any $j = 0, \ldots, i$, $u_j=v_{j-k}\sigma_{j-k}\ldots\sigma_{j-1}v_{j}$ is a node of $\Deb$ and for any $j = 0, \ldots, i-1$, there is an edge $(u_j, u_{j+1}, \sigma_{j})$ in $\Deb$.
    Therefore $\sigma_{0} \ldots \sigma_{i-1} \in \DPaths[i]$.

    If $\sigma_0 \ldots \sigma_{i-1} \in \DPaths[i]$, there exists nodes $u_0, u_1, \ldots, u_i$ of $\Deb$ such that
    $u_0, \sigma_0, u_1, \sigma_1, \ldots, \sigma_{i-1}, u_i$ is a path of $\Deb$.
    Let $v_{-k}, \ldots, v_i$ be the nodes of $\Autom$ and $\sigma_{-k}, \ldots, \sigma_{-1}$ be the symbols such that for any $j = 0, \ldots, i$, $u_j=v_{j-k}\sigma_{j-k}\ldots\sigma_{j-1}v_{j}$.
    By \cref{def:debruijn_graph}, $v_0, \sigma_0, v_1, \sigma_1, \ldots, \sigma_{i-1}, v_{i}$ is a path of $\Autom$ hence
    $\sigma_{0} \ldots \sigma_{i-1} \in \Paths[i]$.

%    Now suppose $\Autom = \Deb[k-1]$, and let us show that $\DLang \subseteq \DLang[k-1]$.
%    Take a sequence $\sigma_{0}\ldots\sigma_{i-1} \in \DPaths[i]$.
%    If $i<k$, then the sequence is part of some node $\sigma_{0}\ldots\sigma_{i-1}\ldots\sigma_{k-1}$ in $\Deb$, and by \cref{def:debruijn_graph}, this sequence is accepted by $\DLang[k-1]$, thus proving our result.
%    Suppose $i \geq k$, so the sequence has the form
%	$\sigma_{0}\ldots \sigma_{k-1} \ldots \sigma_{i-1}$.
%	In order to prove by contradiction, assume that the sequence is not accepted by $\Deb[k-1]$.
%	This means that there is some $j$ with, $k \leq j \leq i$ such that
%	$\sigma_{j-k}\ldots\sigma_{j-2} \in \DLang[k-1]$ and
%	$\sigma_{j-k}\ldots\sigma_{j-2}\sigma_{j-1} \not\in \DLang[k-1]$.
%	This implies that node $\sigma_{j-k}\ldots\sigma_{j-2}\sigma_{j-1}$ does not exist in $\Deb$, which is absurd since the sequence
%	$\sigma_{j-k}\ldots\sigma_{j-1} \in \DLang[k]$.
%    \begin{comment}
%	    The sequence has the form
%	    $\sigma_{0}\ldots\sigma_{k-2}\sigma_{k-1} \ldots \sigma_{i}$.
%	    For any $j=k,\ldots,i(-1?)$,
%	    the word $\sigma_{j-k}\ldots\sigma_{j-1}$ is a node in $\Deb$.
%	    This means that the word $\sigma_{j-k}\ldots\sigma_{j-1} \in \DPaths[k-1]$.
%	    Which means that $\Deb[k-1]$ contains the edge
%	    $(\sigma_{j-k}\ldots\sigma_{j-2}, \sigma_{j-k+1}\ldots\sigma_{j-1}, \sigma_{j-1})$.
%	\end{comment}
  \end{myproof}
\end{proposition}

\begin{theorem}\label[theorem]{theo:remove_1_cycle}
  Consider \cref{alg:stabilization} with input $\Mats, \Deb_0$ (resp. $\Mats, \Deb[k+1]_0$) where $\Deb_0$ (resp. $\Deb[k+1]_0$) is the lift of degree $k$ (resp. $k+1$) of a given automaton $\Autom$.
  If $\Oracle_\epsilon(\Mats, \Deb[k]_0)$ and $\Oracle_\epsilon(\Mats, \Deb[k+1]_0)$ are cycles corresponding to the same word,
  then $\LangAutom[h](\Deb_1) \leq \LangAutom[h](\Deb[k+1]_1)$.
%  let ${\Deb}_1, {\Deb[k+1]}_1$ be corresponding \DBAs{} obtained after removing the same cycle $c$ as in step (\ref{alg:stabilization:s2}) of \cref{alg:stabilization}.
%  \ifreport
%    By the properties of the De Bruijn graph, both $\Deb$ and $\Deb[k+1]$ recognize the same language.
%  \fi
%  Then, $h({\DLang}') \leq h({\DLang[k+1]}')$.
%  Additionally, let ${\DLang}'$ denote the language accepted by ${\Deb}'$ and ${\DLang[k+1]}'$ the language accepted by ${\Deb[k+1]}'$.
%  \rj{theres' a problem here: the automata are not uniquely defined: what if two edges can be removed?}

  \begin{myproof}
    Let $e$ be the edge such that $\Deb_1 = \Deb_0 - e$, that is, the edge removed by the algorithm for $\Deb_0$.
    Let $\sigma_{1}\sigma_{2}\ldots\sigma_{k}\sigma_{k+1}\sigma_{k+2}$ be a sub-word of the repetition of the cycle $c$
    and $v_1, v_2, \ldots, v_{k+3}$ be such that
    $$e = (v_1\sigma_{1}v_2\sigma_{2} \ldots \sigma_{k}v_{k+1}, v_2\sigma_{2} \ldots \sigma_{k} v_{k+1}\sigma_{k+1}v_{k+2}, \sigma_{k+1})$$
    and $(v_{k+2}, v_{k+3}, \sigma_{k+2})$ is an edge of $\Autom$.
    % which means that any word with a sub-word $\sigma_{1} \sigma_{2} \ldots \sigma_{k} \sigma_{k+1}$ will be no longer accepted by $\Deb$.
    Let ${\Deb[k+1]_0}'$ be the graph obtained by removing the node $v_1\sigma_{1}v_2\sigma_{2} \ldots \sigma_{k+1}v_{k+2}$ in $\Deb[k+1]_0$.
    The two automata ${\Deb_1}$ and ${\Deb[k+1]_0}'$ recognize the same language.
    Let ${\Deb[k+1]_0}''$ be the graph obtained by removing the edge $e' = (v_1\sigma_{1}v_2\sigma_{2} \ldots \sigma_{k+1}v_{k+2}, v_2\sigma_{2}v_3\sigma_{3} \ldots \sigma_{k+2}v_{k+3}, \sigma_{k+2})$ in $\Deb[k+1]_0$.
    The language recognized by ${\Deb[k+1]_0}'$ is a subset of the language recognized by ${\Deb[k+1]_0}''$.

    Moreover, as $e'$ is an edge of the cycle $\Oracle_\epsilon(\Mats, \Deb[k+1]_0)$,
    $\LangAutom[h]({\Deb[k+1]_0}'') \leq \LangAutom[h](\Deb[k+1]_1)$.
    Therefore
    \[
      \LangAutom[h]({\Deb}_0) = \LangAutom[h]({\Deb[k+1]_0}') \leq \LangAutom[h]({\Deb[k+1]_0}'') \leq \LangAutom[h](\Deb[k+1]_1).
    \]
%    A way of making $\Deb[k+1]$ not accept any word with a sub-word $\sigma_{1} \sigma_{2} \ldots \sigma_{k} \sigma_{k+1}$ is to remove all its edges $(\sigma_{1} \sigma_{2} \ldots \sigma_{k} \sigma_{k+1}, \sigma_{2} \ldots \sigma_{k} \sigma_{k+1} \sigma_{k+2}, \sigma_{k+2})$, for any $\sigma_{k+2} \in \Mats$.
%    Since, before removing any edges, $\Deb$ and $\Deb[k+1]$ accepted the same language, and after removing the edge they reject the same subset of the language, $\Deb$ and $\Deb[k+1]$ accept the same language after removing the edges.
%    Furthermore, $\rho(\DAdj) = \rho(\DAdj[k+1])$ holds after removing the edges.
%
%    Since it may not necessarily be the case that $\Deb[k+1]$ has to forbid all words with a sub-word $\sigma_{1} \sigma_{2} \ldots \sigma_{k} \sigma_{k+1}$ in order for the resulting graph to be stable, we have that, for the case of a single edge being removed in $\Deb$, the language accepted by $\Deb$ is contained in the language accepted by $\Deb[k+1]$.
%    Additionally, $\rho(\DAdj) \leq \rho(\DAdj[k+1])$.
  \end{myproof}
\end{theorem}

\ifreport
	\claudio{The results commented out below attempt to prove stuff about the edges chosen across different degrees of the DeBruijn graph. Maybe for the journal version?}

\fi

In \cref{sec:application} we show results corroborating \cref{theo:remove_1_cycle}.

\begin{comment}
  \subsection{Notes}
  
  Why don't we backtrack? Why should we?
  
  Given an automaton which represents the decisions that can be made by an orchestration algorithm, we build a \DBA{} of it.
  We start with $n=2$.
  Then we give the graph to the JSR box, which computes the upper and lower bounds on the JSR.
  If the JSR box finds an unstable sequence of matrices, it will return it, and we can identify the path in the graph that represents that sequence.
  If the sequence is of length smaller than $n$, then we cut an edge of that sequence.
  
  In order to pick the right edge, we use, as measure, the grow rate of the JSR of the graph without that edge.
  Alternatively, we can explore the use of some domain specific measure (e.g., in the case of co-simulation, we can remove the edge that has the biggest impact in the performance).
  
  If the sequence is larger than, or equal to, $n$, then we lift the graph.
  
  The larger $n$, the easier it is to prove its stability, and the less conservative the algorithm is.
  
  Can we show that, if the iteration algorithm terminates with an automaton which does not accept infinite paths, then the system is not stabilizable?
\end{comment}

% \ifreport
% 	\input{consensus}
% \fi

\ifreport
	\section{Implementation Details \& Optimality}
\label{sec:impl_opt}

\subsection{Implementation}

The implementation of the stabilization of a CSS is summarized as follows:
\begin{compactenum}
\item find all unstable cycles of length up to 3 using brute force enumeration;
\item since several cycles can be disallowed by removing a single edge, select the edge that disallows the largest number of unstable cycles, and use the entropy of the resulting graph to break ties;
\item repeat steps 1--2 until all allowed cycles have a spectral radius below 1;
\item Use the method of \cite{Legat2017} to determine whether the resulting CSS is stable or whether there is an unstable cycle.
\item if there is an unstable cycle, select the edge that maximizes the entropy of the resulting system  (steps 1--2 of \cref{alg:stabilization});
\item repeat steps 4--5 until the resulting system is stable.
\end{compactenum}

It is easy to see that this implementation is a realization of \cref{alg:stabilization}.
Steps 1--4 are an optimization since they execute relatively quickly, and make the execution of the method in \cite{Legat2017} take less time.

In Step 6, instead of computing the entropy, we compute the spectral radius of the adjacency matrix of the resulting system.
This is equivalent to maximizing the entropy (see \cref{sec:computation_entropy}).

\subsection{Optimality}
\label{sec:optimality}

The solution attained by \cref{alg:stabilization} is not necessarily the optimal solution.
For once, applying different lift degrees will yield different optimal solutions.
Second, \cref{alg:stabilization} removes an edge before finding the next unstable cycle, which means that it misses the chance of optimizing which edge to remove, when more cycles are available (recall steps 1--3 of the above implementation).

Unfortunately, we found no way of guessing which lift degree yields the optimal solution.
However, with small enough constrained switched systems (in the number of matrices and states), it is possible to find the optimal solution, for a given lift degree $k$.

To find the optimal solution, suppose that, for a CSS with a lift degree $k$, we know what the unstable cycles are.
Now we can iterate over all possible procedures to disallow these cycles in the CSS (each procedure is a sequence of edges to be removed), and compute the entropy of the resulting CSS.
The optimal solution is the one that has the maximal entropy.

In order to collect all the unstable cycles, the following procedure can be used:
\begin{compactenum}
\item given a CSS with a lift of degree $k$, apply \cref{alg:stabilization} to find an admissible language, and record all the cycles that were removed throughout the procedure;
\item iterate over all possible ways of disallowing the cycles on the original CSS with a lift of degree $k$, and apply the one that results in a language with maximal entropy;
\item the resulting language is not necessarily admissible, because the best procedure is not necessarily that same as the one picked by \cref{alg:stabilization} in Step 1, so apply \cref{alg:stabilization} to identify and disallow the remaining cycles, adding these to the set of unstable cycles.
\item now repeat Steps 2--3, collecting more and more unstable cycles, until the resulting language is admissible.
\end{compactenum}

The resulting set of unstable cycles represent all possible unstable cycles, and the admissible language found is the optimal solution.

An example application is described in \cref{sec:cosim}.

\begin{comment}
  Suppose that once an admissible language has been found by \cref{alg:stabilization},
  we collect all the unstable cycles that have been found.
  We can iterate over all the possible ways to disallow these cycles in $\Deb$
  and apply the one that leaves a language with maximal entropy.
  The resulting language is not necessarily admissible but we can apply \cref{alg:stabilization} to it to stabilize it.
  Once it is admissible, we can combine the new unstable cycles found with the ones previously found and do the same procedure again.
  Once this terminates, we are certain to have found the language representable by a lift of degree $k$ with maximal entropy.
\end{comment}

\fi

\section{Application}
\label{sec:application}

We first sketch the application of our algorithm to simulation, and then detail the treatment to co-simulation in \cref{sec:cosim}.
\ifreport
	See \cite{Jungers2009,shorten2007stability} for other applications of switched systems, where our technique could be applied as well.
\fi

\subsection{Simulation}
\label{sec:sim}

Consider the problem of approximating the solution $x(t)$ of the system,
\begin{equation}\label{eq:example_linear_ode}
\dert{x}(t) = \bar{A}  x(t) \text{, with } x(0) = x_0,
\end{equation}
using an adaptive simulation algorithm.
\ifreport
	These methods are useful in situations where, e.g., the error tolerance, or runtime performance, can vary as a function of $\bar{x}(t)$ and $t$ \cite{Cellier2006,Shampine1997,Geart1974}.
	In practice, multi-step variable order methods \cite[Section 4]{Cellier2006} are the most commonly used, but for illustrative purposes, we show a single step method. The same analysis can be done for an variable order multi-step method by converting it to a representation in the form of \sys{eq:switched_system}.
\fi

\begin{example}\label{ex:adaptive_numerical_method}
The approximation $\tilde{x}(t)$ of the solution to \sys{eq:example_linear_ode}, computed by a simulation algorithm that uses different step sizes, and different numerical methods, can be modelled as an unconstrained switched system (\sys{eq:switched_system}), with
\begin{align*}
\Mats &= \set{\tilde{A}_{\mathit{fe}, h}, \tilde{A}_{\mathit{md}, h}, \tilde{A}_{\mathit{rg}, h} \middle| h \in \set{0.001, 0.002}} \\
\tilde{A}_{\mathit{fe}, h} &\triangleq \mathbf{I} + \bar{A}  h \\
\tilde{A}_{\mathit{md}, h} &\triangleq \mathbf{I} + \bar{A}  h + (\bar{A}  h)^2/2 \\
\tilde{A}_{\mathit{rg}, h} &\triangleq  \mathbf{I} + \bar{A}  h + (\bar{A}  h)^2/12 + (\bar{A}  h)^3/6 + (\bar{A}  h)^4/24.
\end{align*}
In the above set $\mathcal{A}$, the matrices correspond respectively to the Forward Euler method, the Midpoint method and the Runge-Kutta method.
\end{example}

\begin{comment}
	Let
	$$\Mats^t = \set{A_{\sigma(t-1)}  A_{\sigma(t-2)} \cdots A_{\sigma(0)} : A_{\sigma(j)} \in \Mats, j = 0, \ldots, t}$$
	be the set of all products induced by switching signals with length $t$.
	Note that $x_{t+1} = A x_{0} : \ A \in \Mats^t$.
\end{comment}

In \cref{ex:adaptive_numerical_method}, if one assumes that \sys{eq:example_linear_ode} is stable, that is,
\ifreport
  $$\lim_{t \to \infty} \norm{x(t)} = 0 \text{ for any } x(0),$$
\else
  $\lim_{t \to \infty} \norm{x(t)} = 0 \text{ for any } x(0),$
\fi
then we need to ensure that the error made by discrete approximation $\tilde{x}(t)$ is dissipated.
For this purpose, we introduce the \emph{error switched system}, whose state variable is $e_t = \tilde{x}(t) - x(t)$, and which can be written as
\begin{aligneq}\label{eq:error_dynamics}
e_{t+h} &= A_{\sigma(t)} e_{t} + L_{\sigma(t)}, \text{ with } \\
L_{\sigma(t)} &\triangleq \pargroup{A_{\sigma(t)} - \exp(\bar{A}h})x(t),
\end{aligneq}
where $L_{\sigma(t)}$ is the local error corresponding to $A_{\sigma(t)} \in \Mats$, and $\Mats$ is defined in \cref{ex:adaptive_numerical_method}.
Neglecting $L_{\sigma(t)}$ yields a switched system, the stability of which determines the dissipation of error.

\ifreport
	\begin{remark}
	To derive \cref{eq:error_dynamics},
	\begin{aligneq}
	e_{t+1} &= \tilde{x}_{t+1} - x(t+h) \\
			&= A_{\sigma(t)}\tilde{x}_t - \exp(\bar{A}h)x(t) \\
			&= A_{\sigma(t)}(e_{t} + x(t)) - \exp(\bar{A}h)x(t) \\
			&= A_{\sigma(t)}e_{t} + (A_{\sigma(t)} - \exp(\bar{A}h))x(t)
	\end{aligneq}
	where $\exp(h\bar{A}) = \sum_{k=0}^\infty \frac{h^k}{k!}\bar{A}^k$ is the matrix exponential.
	\end{remark}
\fi

\begin{comment}
\begin{aligneq}
e_t &= \bar{x}_t - x_t \\
e_t + x_t &= \bar{x}_t  \\
x_t &= \bar{x}_t - e_t  \\
\bar{x}_t &= e^{At}  \bar{x}_0 \\
\bar{x}_{t+1} &= e^{A(h+t)}  \bar{x}_0 = e^{Ah}  e^{At}  \bar{x}_0 \\
 			  &= e^{Ah}  \bar{x}_t \\
x_{t+1} &= \tilde{A} x_t \\
e_{t+1} &= \bar{x}_{t+1} - \tilde{A} x_t \\
 		&= \bar{x}_{t+1} - \tilde{A} \bar{x}_t + \tilde{A} e_t \\
 		&= \tilde{A} e_t + e^{Ah} \bar{x}_t - \tilde{A} \bar{x}_t  \\
 		&= \tilde{A} e_t + e^{Ah} \bar{x}_t - \tilde{A} \bar{x}_t  \\
		& \text{target:} \\
e_{t+1} &= \tilde{A} e_t + \text{Local Error} \\
\text{Local Error} &= e^{Ah}\bar{x}_t - \tilde{A} \bar{x}_t
\end{aligneq}

\begin{aligneq}
\lim_{h \to 0} e^{hA} - \tilde{A} = \mathbf{I} - \mathbf{I} - \mathbf{0} = \mathbf{0}
\end{aligneq}
\end{comment}

%\sys{eq:switched_system} is stable when, for any switching signal $\sigma$, and any $x_0 \in \mathbb{R}^n$,
%$
%\lim_{t \to \infty} \norm{ A_{\sigma(t)} \cdots A_{\sigma(h)} A_{\sigma(0)} x_{0}} = 0.
%$

As an example, the error of the adaptive simulation algorithm introduced in \cref{ex:adaptive_numerical_method} may not be stable for a switching signal $111\ldots$, but may be stable for  $2121\ldots$
As a result, the adaptive simulation method may opt for a policy sequence that requires fewer model evaluations (compared to a non-adaptative algorithm), whilst preserving the stability.
\ifreport
	Another way of stating this is to observe that the spectral radius of $\tilde{A}_1$ is larger than 1, that is, $\rho(\tilde{A}_1) > 1$.
	This \emph{does not} imply that the product of any switching signal of the form $2121\ldots$ causes the system to be unstable.

	\begin{figure}[tbh]
	\begin{center}
	  \includegraphics[width=0.8\textwidth]{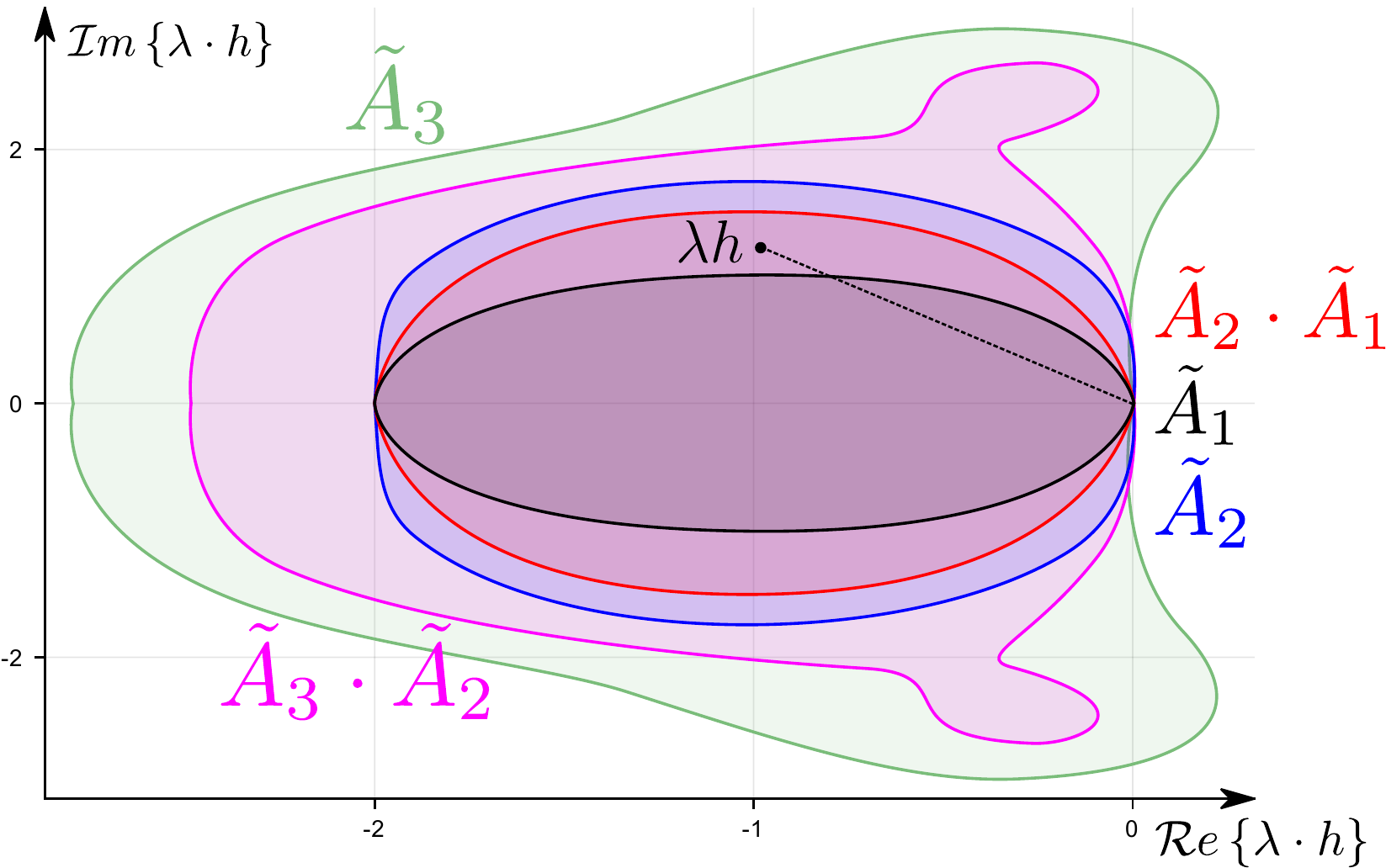}
	  \caption{Domain of numerical stability for some hybrid methods in $\Mats^2$, defined in \cref{ex:adaptive_numerical_method}.}
	  \label{fig:stability_plots}
	\end{center}
	\end{figure}

	To illustrate this fact, \cref{fig:stability_plots} shows the domain of numerical stability for multiple ``hybrid'' solvers applied to the linear system in \cref{eq:example_linear_ode}.
	Each solver is constructed by combining two matrices in $\Mats$, where $\Mats$ is defined in \cref{ex:adaptive_numerical_method}.
	As the figure shows, it is possible that, for a given $h$ and Eigenvalue $\lambda$ of $\bar{A}$, $\lambda  h$ is located outside the stability domain of $\tilde{A}_1$ (shaded area in black), but inside the stability region of the hybrid method $\tilde{A}_2  \tilde{A}_1$ (shaded in red).
	In that case, forbidding any switching signal of the form $11\ldots$, may ensure that the system introduced in \cref{ex:adaptive_numerical_method} is stable.
	
	See \cite[Section 2.4]{Cellier2006} for an example of how to construct the stability domain.
\fi

\subsection{Co-simulation}
\label{sec:cosim}

We apply our algorithm to certify an adaptive orchestration algorithm for the co-simulation of a controlled inverted pendulum.
We will consider two simulators.

In the context of co-simulation, time is discretized into a countable set $T = \set{t_0, t_1, t_2, \ldots} \subset \setreal$, where $t_{i+1} = t_i + H_i$ is the time at step $i$ and $H_i$ is the communication step size at step $i$, with $i=0,1,\ldots$.
From time $t_i \to t_{i+1}$, the simulator $S_j$, with $j=1,2$, is a mapping,
\begin{aligneq}\label{eq:simulator_dts}
\tilde{x}_j(t_{i+1}) &= F_j(t_i,\tilde{x}_j(t_i),u_j(t_i)) \\
y_j(t_i) 	 &= G_j(t_i,\tilde{x}_j(t_i),u_j(t_i))
\end{aligneq}
with state vector $\tilde{x}_j$, input vector $u_j$ and output vector $y_j$.

Simulators exchange outputs only at times $t \in T$.
We assume without loss of generality that the two simulators are coupled in a feedback loop, that is,
\begin{equation}\label{eq:couplings}
u_1 = y_2 \text{ and } u_2 = y_1.
\end{equation}

In the interval $t \in \brackets{t_i, t_{i+1}}$, each simulator $S_j$ approximates the solution to a linear ODE,
\begin{aligneq}\label{eq:linear_ode}
\dert{x}_j &= A_j  x_j + B_j  u_j \\
y_j 	 &= C_j  x_j + D_j  u_j
\end{aligneq}
where $A_j, B_j, C_j, D_j$ are matrices, and the initial state $x_j(t_i)$ is computed in the most recent co-simulation step.
To avoid algebraic loops and keep the exposition short, we assume that either $D_1$ or $D_2$ is the zero matrix. 
Let $D_2 = \mathbf{0}$. 

Since the simulators only exchange outputs at times $t_i, t_{i+1} \in T$, the input $u_j$ has to be extrapolated in the interval $[t_i, t_{i+1})$.
In the simplest co-simulation strategy\footnote{The derivation presented can be applied to more sophisticated input extrapolation techniques, see \cite[Equation~(9)]{Busch2016}.}, this extrapolation is often implemented as a zero-order hold:
$\tilde{u}_j(t) = u_j(t_i)$, for $t \in [t_i, t_{i+1})$.
Then, \cref{eq:linear_ode}
can be re-written to represent the unforced system being integrated by each simulator:
\begin{aligneq}\label{eq:linear_ode_unforced}
\vectorTwo{\dert{x}_j}{\dert{\tilde{u}}_j} &= \vectorTwo{A_j & B_j}{\mathbf{0} & \mathbf{0}}  \vectorTwo{x_j}{\tilde{u}_j}
\end{aligneq}

We can represent the multiple internal integration steps of \sys{eq:linear_ode_unforced}, performed by the simulator $S_j$ in the interval $t \in \brackets{t_i, t_{i+1}}$, as
\begin{aligneq}\label{eq:solver_internal_steps}
\vectorTwo{\tilde{x}_j(t_{i+1})}{\tilde{u}_j(t_{i+1})} &= \tilde{A}^{k_j}_j  \vectorTwo{\tilde{x}_j(t_i)}{\tilde{u}_j}
\end{aligneq}
where, e.g., $\tilde{A}_j = \mathbf{I} + h_j \vectorTwo{A_j & B_j}{\mathbf{0} & \mathbf{0}}$ for the Forward Euler method,
$k_j=(t_{i+1} - t_i)/h_j$ is the number of internal steps, 
and $0 < h_j \leq H_i$ is the internal fixed step size that divides $H_i$.
\ifreport
	Note that this equation implements the mapping in \cref{eq:simulator_dts}.
\fi

At the beginning of the co-simulation step $i$, $u_1(t_i) = y_2(t_i)$ and $u_2(t_i) = y_1(t_i)$, as in \cref{eq:couplings}.
This, together with \cref{eq:linear_ode}, gives,
\begin{aligneq}\label{eq:io_couplings}
u_1(t_i) &= C_2  \tilde{x}_2(t_i) \\
u_2(t_i) &= C_1  \tilde{x}_1(t_i) + D_1 C_2  \tilde{x}_2(t_i).\\
\end{aligneq}

Equations \eqref{eq:linear_ode_unforced}, \eqref{eq:solver_internal_steps}, and \eqref{eq:io_couplings} can be used to represent each co-simulation step by a linear mapping
{\smaller
\begin{aligneq}
\notag
\vectorTwo{\tilde{x}_1(t_{i+1})}{\tilde{x}_2(t_{i+1})} &=
	\underbrace{
	\vectorTwo{\mathbf{I} & \mathbf{0} & \mathbf{0} & \mathbf{0}}{\mathbf{0} & \mathbf{0} & \mathbf{I} & \mathbf{0}}
	\vectorTwo{\tilde{A}^{k_1}_1 & \mathbf{0}}{\mathbf{0} & \tilde{A}^{k_2}_2}
  \vectorFour{\mathbf{I} & \mathbf{0}}{\mathbf{0} & C_2}{\mathbf{0} & \mathbf{I}}{C_1 & D_1  C_2}
	}_{\tilde{A}}
	\vectorTwo{\tilde{x}_1(t_{i})}{\tilde{x}_2(t_{i})}
\end{aligneq}
}
whose state transition matrix $\tilde{A}$ depends on the following items \cite{Gomes2017}:
\ifreport
	\begin{compactitem}
	\item coupling approach --- we used the Jacobi coupling, but Gauss-Seidel, Strong coupling, and others, can be used \cite{Gu2001};
	\item simulator input approximation --- we used constant extrapolation, but other techniques can be applied \cite{Busch2016};
	\item internal solver method --- we showed the Forward Euler as example, but other methods are available;
	\item internal simulator step size $h_j$ --- this value affects matrix $\tilde{A}_j$ and the value $k_j$;
	\item communication step size $H_i$ --- affects $k_j$.
	\end{compactitem}
\else
	\begin{inparaenum}
	coupling approach, input approximations, numerical methods, step size $h_j$, and communication step size $H_i$.
	\end{inparaenum}
\fi
\ifreport
	As evidenced in \cite{Gomes2017c,Schweizer2015c,Busch2016,Busch2010}, each configuration of these items has different stability properties.
\else
	Each configuration of these items affects the stability differently.
\fi
Moreover, as hinted in \cite{Gomes2017d}, an adaptive algorithm that changes the configuration at runtime, based on varying tolerance and performance requirements, is beneficial \emph{as long as it does not make the co-simulation unstable} (recall \cref{ex:adaptive_numerical_method}).
Therefore, with the present work we generate a stabilized CSS, that encodes the set of all possible sequences of configurations that make the co-simulation stable, which can then be consulted during the co-simulation, with little computational cost \cite{Gomes2017d}.

In order to represent an adaptive co-simulation, let $\Mats$ be a set of co-simulation state transition matrices, each representing a particular configuration of the above items.
An adaptive orchestration algorithm will, at the beginning of each co-simulation step $i$, inspect the state variables, and/or local error estimators \cite{Arnold2014}, and decide which configuration should be used to proceed to step $i+1$.
\ifreport
	If we assume that any configuration can be chosen at each co-simulation step, then we can represent the adaptive orchestration algorithm as a switched system.
	On the other hand, if the choice of a configuration for step $i$ depends (in addition to the run-time information) on the configuration chosen for steps $i-1, i-2, \ldots$, then the adaptive orchestration algorithm can be modelled as a constrained switched system.
\fi

Consider the co-simulation of an inverted pendulum that is kept at the equilibrium point using a state feedback controller.
Simulator $S_1$ represents the controller, and simulator $S_2$ represents the pendulum.

Around the equilibrium point, the pendulum can be approximated as a system of the form of \cref{eq:linear_ode}, with
\begin{aligneq}
\notag
A_2 &= \vectorFour{0    &  1       &       0     &      0}%
				  {0 & -(I+ml^2)(b/p) & (m^2gl^2)/p  & 0}%
				  {0    &  0       &       0      &     1}%
				  {0 & -(mlb)/p    &   mgl(M+m)/p & 0} \\
B_2 &= \vectorOne{ 0& (I+ml^2)/p & 0 & ml/p }^\Tr \\
C_2 &= \mathbf{I} \hspace{2em} D_2 = \mathbf{0}
\end{aligneq}
and parameters
$M = 0.5,
m = 0.2,
b = 0.1,
I = 0.006,
g = 9.8,
l = 0.3$.

The controller is a linear quadratic regulator, which, put in the form of \cref{eq:linear_ode}, is
\begin{align}
\notag
A_1 &= \mathbf{0} \hspace{2em} B_1 = \mathbf{0} \\
\notag
C_1 &= \mathbf{0} \hspace{2em} D_1 = \vectorOne{1.0000 & 1.6567 & -18.6854 & -3.4594}.
\end{align}

Assume we can use the Forward Euler and Midpoint methods, with internal fixed step sizes in the set $\set{0.01, 0.02, 0.1, 0.2}$.
\ifreport
	Furthermore, the co-simulation step size can be $H=0.1$ or $H=0.2$.
\else
	Furthermore,  $H=0.1$ or $H=0.2$.
\fi
\ifreport
	Note that the internal step sizes must always divide the communication step size $H$, and that the numerical method and step size used in the controller simulator have no impact in the co-simulation stability, because it has no internal dynamics.
\fi
Then, applying Equations (\ref{eq:linear_ode_unforced}), (\ref{eq:solver_internal_steps}), and (\ref{eq:io_couplings}), we get a switched system over 8 matrices.

The matrices $A_2$ (corresponding to $H=0.2$, $h_1=0.2$, Forward Euler),
$A_3$ (corresponding to $H=0.2$, $h_1=0.02$, Midpoint) and
$A_4$ (corresponding to $H=0.2$, $h_1=0.02$, Forward Euler)
have a spectral radius larger than one.
This means that the switching signals $222...$, $333...$ and $444...$ should be forbidden.

Applying \cref{alg:stabilization} directly to the unconstrained switched system 
\ifreport
	(which corresponds to a lift of degree 0)%
\fi
, leads to removal of the edges with labels $2, 3$ and $4$.
This completely disallows the use of the matrices $A_2, A_3$ and $A_4$.
The resulting language turns out to be admissible, its entropy is $\log_2(5)$.

Applying \cref{alg:stabilization} to a lift of degree 1, we get a constrained switched system with the automaton shown in \cref{fig:adaptive_cosim_stable}, where the edges in red were removed by the algorithm.
We can see that the matrices $A_2$, $A_3$ and $A_4$
are now allowed by the algorithm (only the cyclic application of each one of these matrices is still disallowed).
This solution is less conservative than the one with degree 0.
\ifreport
	Its entropy is $~log_2 (7.26)$.
\fi
One allowed cycle is $32645$, where the symbols 5 and 6 seem to play the role of stabilizing the cycle.

\ifreport
  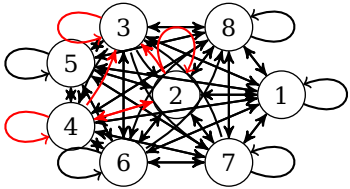
\begin{figure}[htb]
    \centering
    \begin{tikzpicture}[xscale=0.7,yscale=0.6]
  \Vertex[L={1},x=2,y=0]{1};
  \Vertex[L={8},x=1,y=1.5]{8};
  \Vertex[L={3},x=-1,y=1.5]{3};
  \Vertex[L={5},x=-2,y=0.7]{5};
  \Vertex[L={4},x=-2,y=-0.7]{4};
  \Vertex[L={6},x=-1,y=-1.5]{6};
  \Vertex[L={7},x=1,y=-1.5]{7};
  \Vertex[L={2},x=0,y=0]{2};
  \SetUpEdge[style={post,{->}, thick}]
  \Edge(1)(2)
  \Edge(1)(3)
  \Edge(1)(4)
  \Edge(1)(5)
  \Edge(1)(6)
  \Edge(1)(7)
  \Edge(1)(8)
  \Edge(2)(1)
  \Edge(2)(5)
  \Edge(2)(6)
  \Edge(2)(7)
  \Edge(2)(8)
  \Edge(3)(1)
  \Edge(3)(5)
  \Edge(3)(6)
  \Edge(3)(8)
  \Edge(4)(1)
  \Edge(4)(5)
  \Edge(4)(6)
  \Edge(4)(7)
  \Edge(4)(8)
  \Edge(5)(1)
  \Edge(5)(2)
  \Edge(5)(3)
  \Edge(5)(4)
  \Edge(5)(6)
  \Edge(5)(7)
  \Edge(5)(8)
  \Edge(6)(1)
  \Edge(6)(2)
  \Edge(6)(3)
  \Edge(6)(4)
  \Edge(6)(5)
  \Edge(6)(7)
  \Edge(7)(1)
  \Edge(7)(2)
  \Edge(7)(4)
  \Edge(7)(5)
  \Edge(7)(6)
  \Edge(7)(8)
  \Edge(8)(1)
  \Edge(8)(2)
  \Edge(8)(3)
  \Edge(8)(4)
  \Edge(8)(5)
  \Edge(8)(7)
  \SetUpEdge[style={post,bend right=20,thick}]
  \Edge(3)(7)
  \Edge(7)(3)
  \Edge(6)(8)
  \Edge(8)(6)
  \SetUpEdge[style={post,thick,color=red}]
  \Edge(2)(4)
  \Edge(4)(2)
  \SetUpEdge[style={post,bend right=8,thick}]
  \Edge(3)(2)
  \Edge(3)(4)
  \SetUpEdge[style={post,bend right=8,thick,color=red}]
  \Edge(2)(3)
  \Edge(4)(3)
  \draw[thick,->] (1) to[out=-25,in=25,looseness=7] (1);
  \draw[thick,->,color=red] (2) to[out=115,in=65,looseness=10] (2);
  \draw[thick,->,color=red] (3) to[out=155,in=205,looseness=7] (3);
  \draw[thick,->,color=red] (4) to[out=155,in=205,looseness=7] (4);
  \draw[thick,->] (5) to[out=155,in=205,looseness=7] (5);
  \draw[thick,->] (6) to[out=155,in=205,looseness=7] (6);
  \draw[thick,->] (7) to[out=-25,in=25,looseness=7] (7);
  \draw[thick,->] (8) to[out=-25,in=25,looseness=7] (8);
  \SetUpEdge[color=red]
\end{tikzpicture}
    \caption{Solution with entropy $\log_2(7.2568898)$.}
    \label{fig:adaptive_cosim_stable}
  \end{figure}
\fi

We applied \cref{alg:stabilization} to the lifts of degree 0, 1 and 2. % of the switched system of increasing degrees.
At each application of the algorithm, a stable constrained switched system was produced, with an entropy that increased with the degree of the lift.
These results, summarized in \cref{tab:entropy_over_k}, corroborate \cref{theo:remove_1_cycle}.

\begin{table}[tbh]
  \centering
  \caption{Entropy achieved per lift degree.}
  \label{tab:entropy_over_k}
  \begin{tabular}{l|l|l}
    $k$ & Entropy [bit] & CPU time [s]\\
    \hline
    0 & $\log_2(5)$ & 0.13\\
    1 & $\log_2(7.2568898)$ & 1.8\\
    2 & $\log_2(7.7083039)$ & 280
  \end{tabular}
\end{table}

\ifreport
  This application to co-simulation illustrates an important advantage of the method presented in \cite{Legat2017}: it is capable of finding large unstable cycles.
  This method does not find the unstable cycles by iterating through the cycles of some length $K$ but instead extracts them from an infinite switching signal,
  hence it is not harder for the method to find large unstable cycles.
  For the lift of degree 2 for example, it found the unstable cycle $542245332$ of length 9, and in a subsequent iteration found the cycle $224533542245335422453354224523254$ of length 33.
  A brute force method would have to enumerate all $8^{33}$ cycles to achieve the stabilization of the adaptive solver.
  
  Regarding the optimality of the solution found for the lift with degree 1, we have applied the procedure detailed in \cref{sec:optimality} to confirm that $\log_2(7.2568898)$ is indeed the maximal entropy for that degree.

\else
  In the implementation,
  we first find all unstable cycles of length from 1 to 3 by iterating over all cycles of these lengths using brute force and selecting the ones that have a spectral radius above one.
  Note that several cycles can be disallowed at the same time using the same edge.
  Therefore we remove in priority the edges breaking the most cycles and use the entropy to break ties.
  We then use the method of \cite{Legat2017} to determine whether the system is stable or whether it can provide an unstable cycle.
  
  This application to co-simulation illustrates an important advantage of the method presented in \cite{Legat2017}: it is capable of finding large unstable cycles.
  For the lift of degree 2 for example, it found the cycle $224533542245335422453354224523254$ of length 33.
\fi

\begin{comment}
  \subsection{Notes}
  
  \claudio{Is the set of matrices used by the adaptive simulator commonly reducible?}
  
  For long running co-simulations, it is important that we certify that all decisions that can be made by an orchestration algorithm always produce stable trajectories.
  
  Our approach is to take the existing policies of an orchestrator and restrict them in a way that ensures that no unstable co-simulation can be computed, for any initial condition \claudio{for non linear systems, this is not the case.}.
  
  I think we can use the pendulum control system as an example.
  
  \claudio{Can we draw stability plots for co-simulation?}
  
  % Example from here: http://ctms.engin.umich.edu/CTMS/index.php?example=InvertedPendulum&section=SystemModeling
  
  Co-simulation as a switched control problem.
  
  We assume that the origin of the original system is stable.
\end{comment}

\section{Related Work}
\label{sec:related_work}

\ifreport
	There has been a huge effort to understand how to ensure that a (discrete time) switched system is stable.
	We refer the interested reader to \cite{Lin2009,Sun2005,Jungers2009,parrilo2007approximation,ahmadi2011analysis} for introductions and surveys on this subject.
\else
	There has been a huge effort to understand how to ensure that a (discrete time) switched system is stable \cite{Lin2009,Jungers2009,parrilo2007approximation,ahmadi2011analysis}.
\fi
To the best of our knowledge, the problem we introduce here%
\ifreport
	, i.e. finding the largest set of switching signals that guarantees the stability of the system,
\fi
has never been studied.
In the broader field of stabilization of switched systems, we can highlight the works in
\ifreport
  \cite{GuishengZhai2002,Xu2002,Xu2004,Zhang2009a,Hetel2006,Zhao2012,Pettersson2003,kundu2017stabilizing,prabhakar2017formal}%
\else
  \cite{Xu2004,Hetel2006,Zhang2009a,Zhao2012,kundu2017stabilizing,prabhakar2017formal}%
\fi%
.
The key difference with our work is the goal: we are not satisfied with a single stable switching signal; we want to provide
the maximum flexibility to the stabilized CSS, 
which can make use of this flexibility to choose the most appropriate switching signal.
%that would choose the most appropriate stable switching sequence based on the performance constraints.
\ifreport
  The works in \cite{Hetel2006,Xu2002,Xu2004,Pettersson2003,kundu2017stabilizing,prabhakar2017formal} are focused on continuous time systems, and \cite{Zhao2012,Xu2002,Xu2004} aim at deriving state feedback laws (in addition to switching signals) that make the system stable.
  
  The approach followed in \cite{Zhao2012} assumes that each mode of the system is stable.
  In our case, the goal is the same but we tolerate unstable modes.
  
  The approach in \cite{GuishengZhai2002} is interesting because it allows the combination of stable and unstable modes, in order to ensure stability.
  However, no algorithm is provided to find these combinations.
  
  The aim of \cite{kundu2017stabilizing} is different as it describes the search for one particular stable trajectory
  while we maximize the size of a language of stable switching signals.
  
  \cite{Sloth2014} describes the stability analysis for continuous switched systems with parametric uncertainties. 

  \cite{Zhang2009a} focuses on proving that a system is stabilizable, rather than making the system stable. 
  It deals with forced discrete time switched systems, and the stabilization procedure finds a control policy (switching signal, and input) that stabilizes the system.
  This is in contrast to our goal, which is to find all policies that make the system stable, and maximize this set.
\fi
In the context of co-simulation, the reader can find stability analysis of traditional orchestration algorithms in 
\ifreport
	\cite{Schweizer2015d,Schweizer2015c,Busch2016,Gomes2017c}.
\else
	\cite{Busch2016,Gomes2017c} and references thereof.
\fi

\section{Conclusion}
\label{sec:conclusion}

We introduce a new problem in the context of constrained switched
systems: 
1) to restrict the switching possibilities of
the system, so as to ensure its stability, and 2) to leave as many
switching policies as possible (provided that the system becomes stable).

The motivation for leaving as many switching policies as possible lies in the fact that, in adaptive co-simulation, the orchestration algorithm will make the best possible choice as a function of information obtained during the simulation.
\ifreport
	We restrict the switching possibilities to be representable by an automaton because of their great efficiency.
\fi

The problem is interesting in that it transforms a control
problem into the problem of building an optimal language, that is,
optimizing the construction of an automaton. 
By combining classical control concepts for switched systems 
\ifreport
	(like the CJSR)
\fi 
, with classical automata-theoretic concepts 
\ifreport
	(like the entropy of shifts)%
\fi
, one can design algorithms to solve this problem.
Our algorithm takes the form of a hierarchy of sufficient conditions,
where increasingly better solutions are found by lifting the automaton (see \cref{fig:1cycle} and \cref{tab:entropy_over_k}).
Essentially, this allows one to control the optimality of the solution, at the cost of processing power and memory.

This work is aimed to be a proof of concept, and we leave many research questions open.
\ifreport
	We plan to investigate the conservativeness of restricting ourselves to regular languages (see \cref{ex:pumping}). 
	Second, we want to understand how our method can be optimized for the particularizes of co-simulation, and apply it to nonlinear systems.
	Finally, we plan to modify \cref{alg:stabilization} so that stronger theoretical results can be proven.
\else
	We want to investigate the conservativeness of restricting ourselves to regular languages, and to modify \cref{alg:stabilization} so that stronger theoretical results can be proven.
\fi

\ifreport
	\claudio{Future work: how can we find the optimal (something like Benoit's procedure, assuming that we have the maximal length of the cycles in the graph). 
	As a reviewer puts it: ``Is it true that the maximum is achieved when the lifting order k is larger than the maximal length of circles in the graph?''}
	
	\claudio{From the same reviewer: ``The structure of regular language generated by automata is
	simple: prefix + circles + suffice. As it is not hard to
	analyze the spectral radius of each circle, is it possible
	to find the maximal set of switching signals by just
	removing the minimal number of unstable circles and keeping
	the maximal number of stable circles?''}
	
	\claudio{Future work from reviewer: In the only realistic co-simulation example, the
	inverted pendulum, the controller does not have any
	dynamics. What would be the value of using the
	co-simulation? The authors should consider using a more
	realistic example.}
	\benoit{Maybe not add it in the conclusion or rewrite it in a more positive fashion}
	\claudio{Sorry, I placed that here just for future reference (it was not meant to make it to the final version)}
	
\fi

\ifreport
	\section*{Acknowledgments}
	
  This research was partially supported by
  a PhD fellowship grant from the Agency for Innovation by Science and Technology in Flanders (IWT, dossier 151067),  
  the Belgian Interuniversity Attraction Poles, 
  the ARC grant 13/18-054 from Communaut\'e fran\c{c}aise de Belgique, 
  a F.R.S.-FNRS Research Fellowship, and
  Flanders Make vzw, the strategic research centre for the manufacturing industry.
\fi

%\appendix

%\input{notes}

\ifreport
	\bibliographystyle{plain}
	\bibliography{bibliography_benoit,bibliography_claudio}
\else
	\bibliographystyle{IEEEtran}
	\bibliography{rmurl,bibliography_benoit,bibliography_claudio}
\fi

\ifreport
	\appendix
	
	\section{Computation of the Entropy}
\label{sec:computation_entropy}

\subsection{Spectral Radius of Adjacency Matrix}

Consider a given CSS $S = \tuple{\Mats, \Autom}$, and let $B$ be the adjacency matrix of $\Autom$.

The matrix element $b_{ij}$ of $B^k$ gives the number of different paths of length $k$ from node $i$ to node $j$ \cite{West2001}.
Hence, $\bnorm{B^k}$ gives a measure of the size of the matrix comprised by the number of different paths from each node to each other node (see, e.g., \cite[Remark~2]{Legat2016}),
and $\bnorm{B^k}^{\frac{1}{k}}$ gives the growth rate of this quantity.
Taking the limit $k \to \infty$, we have the spectral radius of the adjacency matrix:
$$ \rho(B) = \lim_{k\to\infty} \bnorm{ B^k }^{1/k}. $$

\begin{example}
Recall \cref{ex:edge_removal_example}, let $B_1$ denote the adjacency matrix of the automata in \cref{fig:edge_removal_example} without the edge $v_1 \xrightarrow{2} v_2$, and let $B_2$ denote the adjacency matrix of the same automata, without the edge $v_2 \xrightarrow{3} v_3$.
Then $\rho(B_1) = 1 < \rho(B_2) \approx 1.6180$.
\end{example}

\subsection{Edge Shift}

The logarithm of the spectral radius of the adjacency matrix of an \emph{irreducible} automaton gives the entropy of its \emph{edge shift} \cite[Theorem~4.3.1]{Lind1995}.
An automaton is \emph{irreducible} if for every pair of nodes $u, v$, there exists a path from $u$ to $v$ accepted by the automaton. In other words, the graph consists of a single strong component.

\begin{definition}[{\cite[Definition~2.2.5]{Lind1995}}]
  \label[definition]{def:edgeshift}
  The \emph{edge shift} of an automaton $\Autom = (V, E)$ is the language
  recognized by
    the automaton $\Autom' = (E, E')$ with the transitions
  $((u,v,\sigma),(v,w,\sigma'),(v,w,\sigma')) \in E'$ for each $(u,v,\sigma), (v,w,\sigma') \in E$.
\end{definition}

An edge shift of automaton \cref{fig:even} is illustrated in \cref{fig:edgeshift}.

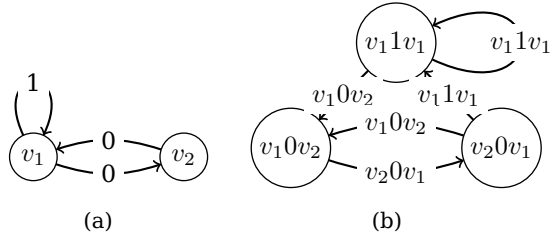
\begin{figure}[!ht]
    \centering
    \begin{subfigure}[t]{0.3\textwidth}
      \centering
      \begin{tikzpicture}
        \Vertex[L={$v_1$}]{1}
        \EA[L=$v_2$,unit=2](1){2}
        \tikzset{EdgeStyle/.style = {->}}
        \tikzset{EdgeStyle/.append style = {bend right=18}}
        \Edge[label=0](1)(2)
        \Edge[label=0](2)(1)
        %\Loop[dist=1cm,dir=EA](3)
        \draw[thick,->] (1) to [out=115,in=65,looseness=10] node [midway, fill=white] {1} (1);
      \end{tikzpicture}
       \caption{ }
       \label{fig:even}
    \end{subfigure}%
    ~
    \begin{subfigure}[t]{0.3\textwidth}
      \centering
      \begin{tikzpicture}
        \Vertex[L={$v_11v_1$},unit=1.4]{1}
        \SOEA[L=$v_20v_1$,unit=1.4](1){2}
        \SOWE[L=$v_10v_2$,unit=1.4](1){3}
        \tikzset{EdgeStyle/.style = {->}}
        \Edge[label=$v_10v_2$](1)(3)
        \Edge[label=$v_11v_1$](2)(1)
        \tikzset{EdgeStyle/.append style = {bend right=18}}
        \Edge[label=$v_10v_2$](2)(3)
        \Edge[label=$v_20v_1$](3)(2)
        \draw[thick,->] (1) to [out=-25,in=25,looseness=10] node [midway, fill=white] {$v_11v_1$} (1);
      \end{tikzpicture}
       \caption{ }
       \label{fig:edgeshift}
    \end{subfigure}
    \caption{An automaton (a) and its edge shift (b).}
  \end{figure}

It turns out that the entropy of the edge shift is equal to the entropy of the language recognized by the automaton if the automaton is \emph{right-resolving}~\cite[Proposition~4.1.13]{Lind1995}.

\begin{definition}[{\cite[Definition~3.3.1]{Lind1995}}]
  An automaton $\Autom$ is \emph{right-resolving} if for every vertex $v$, the outgoing edges have different symbols.
\end{definition}

Every regular language is recognized by a right-resolving automaton.
Moreover, there are automated ways to obtain such an automaton from a starting representation of a language with an automaton that is not right-resolving \cite[Section~3.3]{Lind1995}.

Since the automata considered here are right resolving and irreducible, the entropy is computed by computing the spectral radius of the adjacency matrix of the CSS.

\fi

\end{document}